\documentclass[12pt]{article}
\nonstopmode
\RequirePackage[colorlinks,citecolor=blue,urlcolor=blue,linkcolor=blue]{hyperref}
\hypersetup{
colorlinks = true,
citecolor=blue,
urlcolor=blue,
linkcolor=blue,
pdfauthor = {Alexey Kuznetsov},
pdfkeywords = {Riemann zeta function, Riemann-Siegel formula, Gaussian quadrature, high-precision algorithm},
pdftitle = {Simple and accurate approximations to the Riemann zeta function},
pdfpagemode = UseNone
}
\usepackage{matlab-prettifier}
\usepackage[titletoc,title]{appendix}
\usepackage{graphicx,xspace,colortbl}
\usepackage{amsmath,amsthm,amsfonts,mathrsfs}
\usepackage{color}
\usepackage{fancybox}
\usepackage{epsfig}
\usepackage{subfig}
\usepackage{pdfsync}
    \def\qed{\hfill$\sqcap\kern-8.0pt\hbox{$\sqcup$}$\\}
    \def\beq{\begin{eqnarray}}
    \def\eeq{\end{eqnarray}}
    \def\beqq{\begin{eqnarray*}}
    \def\eeqq{\end{eqnarray*}}

    \def\re{\textnormal {Re}}
    \def\im{\textnormal {Im}}

    \def\r{{\mathbb R}}
    \def\c{{\mathbb C}}
    	
    \def\d{{\textnormal d}}
    \def\i{{\textnormal i}}

    \def\ii{{\mathcal I}}

\newtheorem{proposition}{Proposition}

\theoremstyle{definition}

\title{Simple and accurate approximations to the Riemann zeta function}
\author{ 
{Alexey Kuznetsov\footnote{Dept. of Mathematics and Statistics,  York University,
4700 Keele Street, Toronto, ON, M3J 1P3, Canada.  
  Email: akuznets@yorku.ca}} }

\date{\today}

\begin{document}
\maketitle

\begin{abstract}
We develop approximations to the Riemann zeta function that enable high-precision computation in the critical strip and in other vertical strips. These approximations combine the main sum in the Riemann-Siegel formula with a simple approximation to the remainder term involving only elementary functions and certain precomputed coefficients obtained by Gaussian quadrature. We also derive approximations for the derivative of the Riemann zeta function and present extensive numerical evidence demonstrating the accuracy of these approximations.
\end{abstract}
{\vskip 0.15cm}
 \noindent {\it Keywords}: Riemann zeta function, Riemann-Siegel formula, Gaussian quadrature, high-precision \\ algorithm
{\vskip 0.25cm}
 \noindent {\it 2020 Mathematics Subject Classification}: Primary 11M06, Secondary 11Y35

\section{Introduction and main results}\label{section:Intro}

There exist many methods for computing  the Riemann zeta function. One of the simplest is the Euler-Maclaurin summation method
\cite{Borwein_2000,Johansson_2015,Rubinstein_2005}, which is easy to implement,  provides rigorous error bounds, and allows for high-precision computation of $\zeta(s)$. However, a major drawback of this method is that it requires summing  $O(|s|)$ terms, making it impractical for large values of $s$. Another highly efficient algorithm, developed by Borwein \cite{Borwein_2000b}, also has rigorous error bounds and enables high-precision computation of $\zeta(s)$ when $\im(s)$ is not large.

For values of $s$ with large $\im(s)$, the preferred method is the Riemann-Siegel formula and its various extensions. The original Riemann-Siegel formula \cite{Borwein_2000, Gabcke1979,Rubinstein_2005,Titchmarsh1987} was developed for computing $\zeta(s)$ on the critical line $\re(s)=1/2$. It consists of a main sum of $N_t:=\lfloor \sqrt{t/(2\pi)} \rfloor$ terms (here and throughout this paper we denote $s=\sigma+ \i t$ and assume that $t>0$) and a remainder term, which can  be expressed in terms of certain integrals or expanded as an asymptotic series. Keeping the first $j$ terms of this asymptotic series results in an error of $O(t^{-\frac{1}{4} - \frac{j}{2} })$ (see \cite{Gabcke1979} for rigorous and effective upper bounds for these errors).  Considerable work has been devoted to  extending and improving the Riemann-Siegel formula.  Odlyzko and Sch\"onhage \cite{Odlyzko_1988} introduced a fast algorithm for evaluating $\zeta(s)$ at multiple points, which is particularly useful for locating the zeros of the Riemann zeta function on the critical line. Hiary  
\cite{Hiary_2011} developed an algorithm that reduces the complexity of evaluating $\zeta(1/2+\i t)$ at a single point to $O(t^{4/13+o(1)})$, compared to the standard Riemann-Siegel formula's complexity of $O(t^{1/2})$. A simpler  version of Hiary's algorithm \cite{Hiary_2011} achieves complexity of $O(t^{1/3+o(1)})$ with minimal memory requirements. 
Another algorithm due to Hiary \cite{Hiary_2016} produces Riemann-Siegel-type approximations (requiring $O(t^{1/2+o(1)})$ terms) by starting from the Euler-Maclaurin formula. 
Smoothed versions of the Riemann-Siegel formula have been derived in  \cite{Berry_Keating1992,Kuznetsov2023,Paris_Cang1997,Rubinstein_2005}
and a Riemann-Siegel-type formula for computing $\zeta(s)$ for general values of $s$ (not necessarily on the critical line) was stated in \cite{Odlyzko_1988}, with error estimates developed in \cite{Reyna_2011}.

In this paper, we construct approximations to $\zeta(s)$ that are valid in arbitrary vertical strips, including the critical strip, for large values of $\im(s)$.    Our method is conceptually similar to the approach of Galway \cite{Galway_2001}, whose algorithm computes $\zeta(s)$ by numerically evaluating the integrals that appear in the error term of the Riemann-Siegel formula. Galway's approximations can yield highly  accurate results, but their implementation requires careful tuning, as they depend on three parameters $\{z_1, z_2, M\}$ that must be appropriately chosen to achieve the best accuracy. In contrast, our approximations  depend on a single integer parameter $p$ along with certain precomputed complex coefficients $\{\omega_{p,j}\}_{0\le j \le p}$ and $\{\lambda_{p,j}\}_{1\le j \le p}$. Our approximations allow  high-precision computation of $\zeta(s)$ in the critical strip, and near it, for values of $|\im(s)|>100$. High-precision computations of the Riemann zeta function are not merely of theoretical interest; they can lead to significant mathematical discoveries. For example, computations of the nontrivial zeros of $\zeta(s)$ to high precision were a crucial ingredient in the disproof of Mertens conjecture
\cite{OdlyzkoRiele_1985}. 
\begin{figure}[t]
\centering
\subfloat[][$\lambda_{20,j}$ for $1\le j \le 20$]{\label{fig_1a}\includegraphics[height =5.25cm]{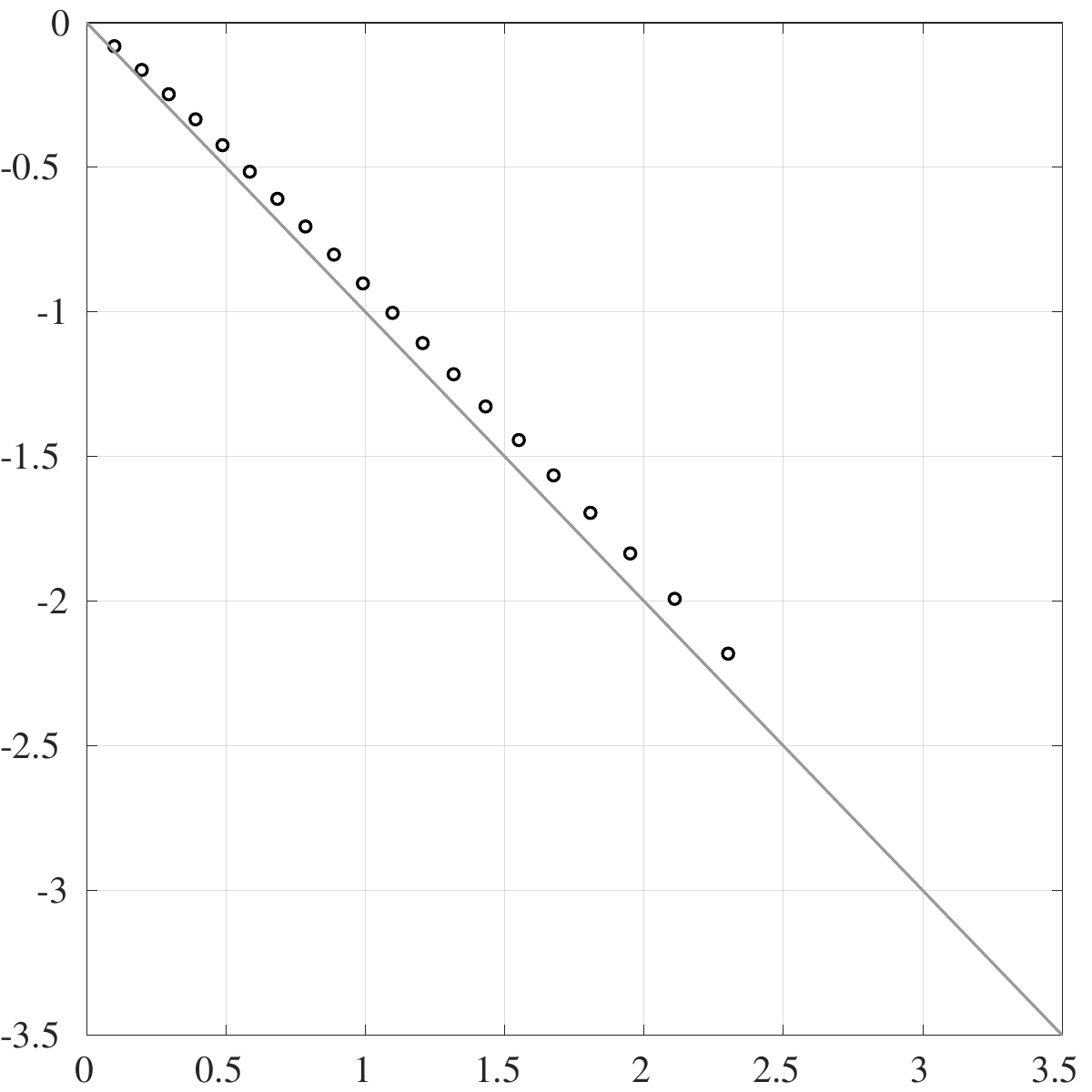}}
\hspace{0.5cm}
\subfloat[][$\lambda_{40,j}$ for $1\le j \le 40$]{\label{fig_1b}\includegraphics[height =5.25cm]{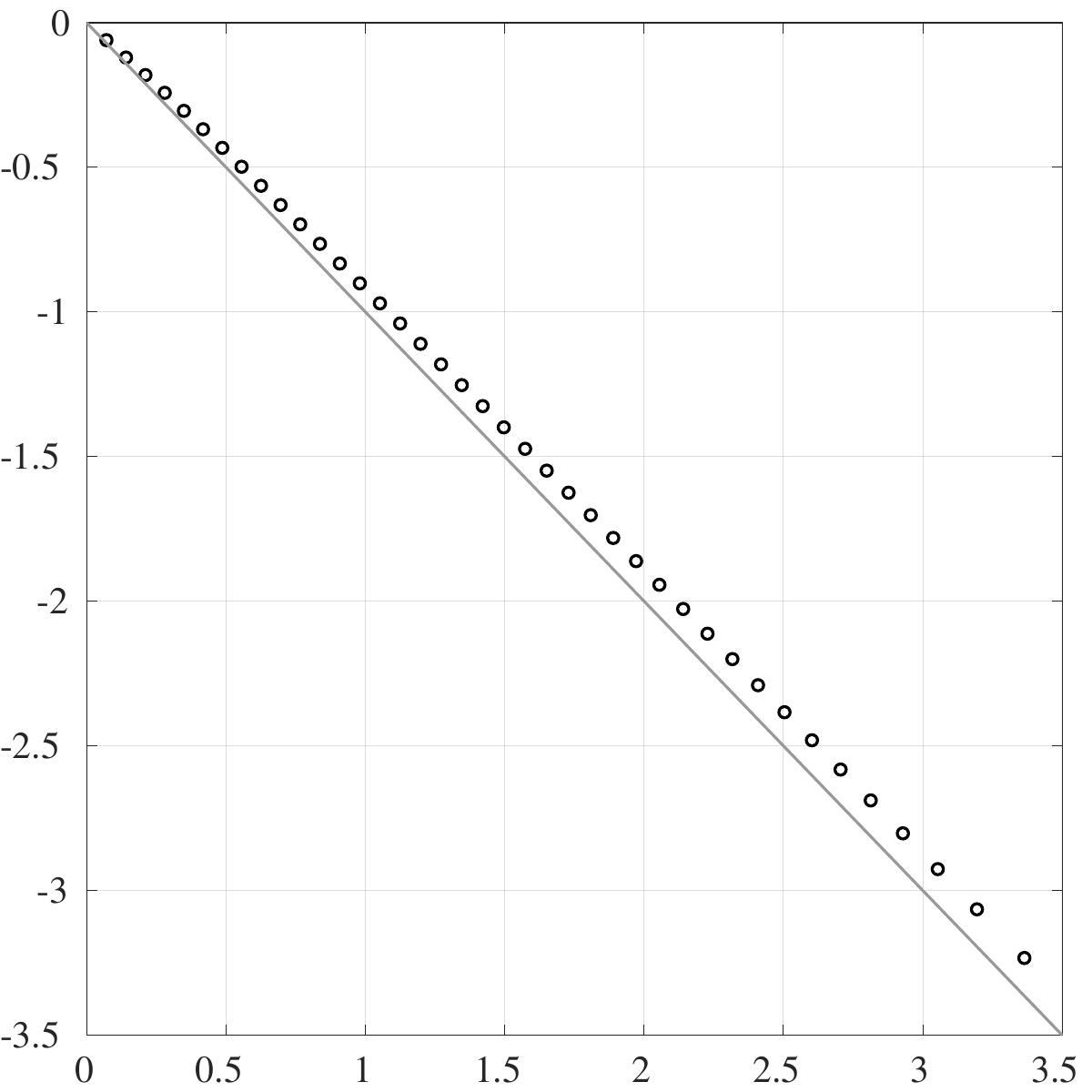}} \\
\subfloat[][$|\omega_{20,j}|$ for $0\le j \le 20$]{\label{fig_1c}\includegraphics[height =5.25cm]{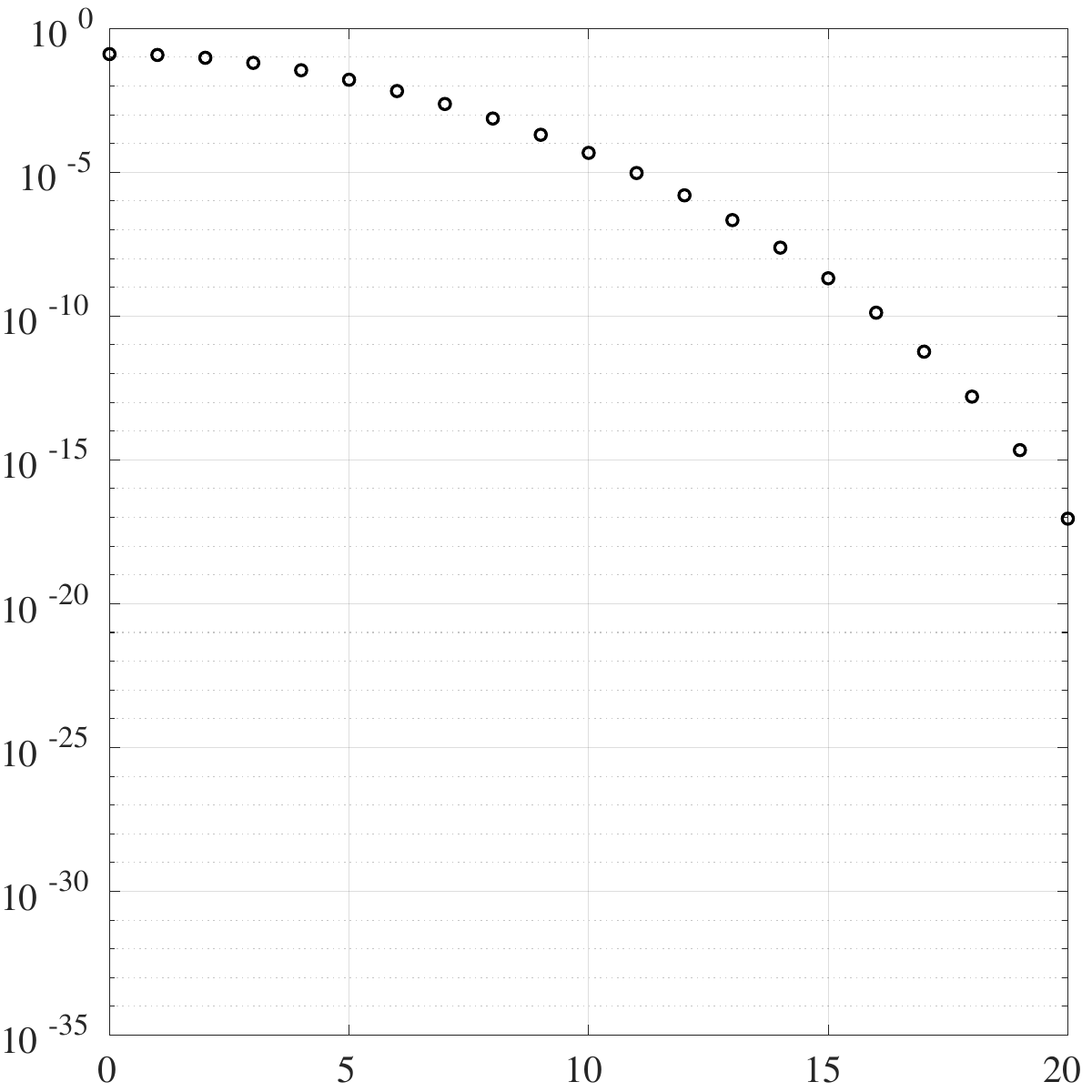}} 
\hspace{0.5cm}
\subfloat[][$|\omega_{40,j}|$ for $0\le j \le 40$]{\label{fig_1d}\includegraphics[height =5.25cm]{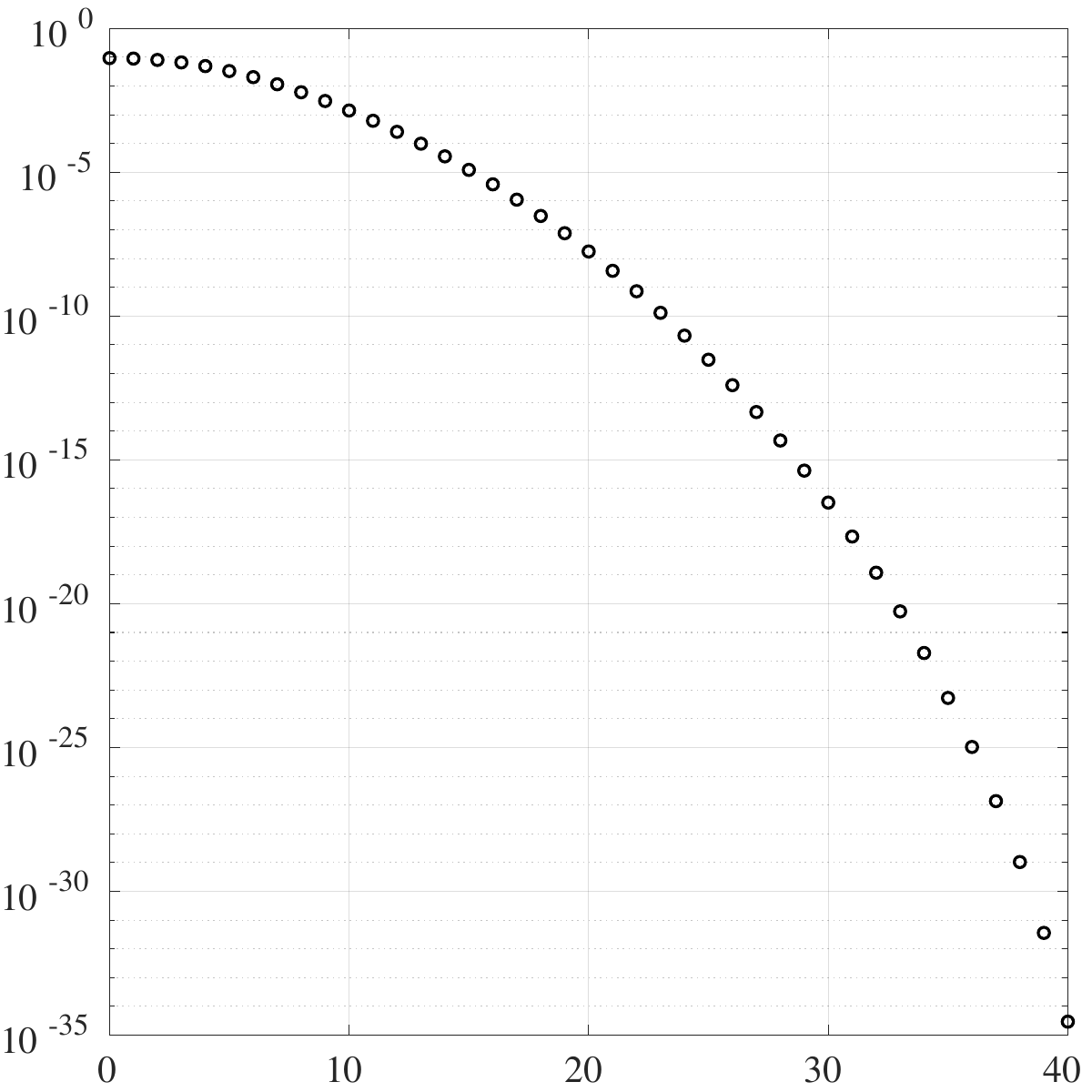}}  
\caption{The graphs (a) and (b) show the complex numbers $\lambda_{p,j}$ for $p\in \{20,40\}$ (all of them are located in the fourth quadrant of the complex plane). The graphs (c) and (d) show the values of $|\omega_{p,j}|$ for $p\in \{20,40\}$ (here $j$ is on the $x$-axis and $y$-axis is in log-scale). } 
\label{fig1}
\end{figure}

The key components of our approximations for the Riemann zeta function are the $2p+1$ complex numbers $\{\omega_{p,j}\}_{0\le j \le p}$ and $\{\lambda_{p,j}\}_{1\le j \le p}$. We explain in Section \ref{section2} how these numbers are derived. We have precomputed $\omega_{p,j}$ and  $\lambda_{p,j}$ to sufficiently high precision for 
 all $1\le p\le 30$ and for many values of $p$ in the range $30<p\le 150$; the results can be downloaded from the author's  \href{https://kuznetsovmath.ca/}{webpage}. 
The values of $\omega_{p,j}$ and $\lambda_{p,j}$ for $p \in \{5,8,10\}$ are listed in Appendices A and B. To illustrate the magnitude and distribution of these numbers in the complex plane, Figure \ref{fig1} displays $\lambda_{p,j}$ and $|\omega_{p,j}|$ for $p\in \{20,40\}$. Note that in Figures \ref{fig_1a} and \ref{fig_1b} we plot the points $\lambda_{p,j}$ in the complex plane (these numbers are indexed in increasing order of absolute value), while in Figures
\ref{fig_1c} and \ref{fig_1d} we plot the real numbers $|\omega_{p,j}|$ against $j$ on the $x$-axis. We observe that the numbers $\lambda_{p,j}$ and $\omega_{p,j}$ are not large in magnitude, that $\lambda_{p,j}$  lie in the fourth quadrant of the complex plane slightly above the ray $\arg(z)=-\pi/4$ and that  $|\omega_{p,j}|$ decay rapidly as $j$ increases. 

Given the numbers $\{\omega_{p,j}\}_{0\le j \le p}$ and $\{\lambda_{p,j}\}_{1\le j \le p}$, we define
\begin{equation}\label{def:ii_Mp}
\ii_{M,p}(s):= \omega_{p,0} M^{-s}   + \sum\limits_{j=1}^p 
\omega_{p,j} \Big[
e^{-2 \pi M   \lambda_{p,j}} 
 \big(M+ \i \lambda_{p,j}\big)^{-s}
+e^{2  \pi M    \lambda_{p,j}} 
 \big(M- \i \lambda_{p,j} \big)^{-s} \Big],
\end{equation}
where  $M>0$ and $s\in \c$.
For a function $f : \c \mapsto \c$  we denote
\begin{equation}\label{def:bar_f}
\bar{f}(s):={\overline {f(\bar s)}}. 
\end{equation}
Clearly, if $f$ is an entire function of $s$, then so is $\bar f$. For $p\in {\mathbb N}$  and $N \in {\mathbb N} \cup \{0\}$ we introduce
\begin{equation}\label{def:F_sNP}
F(s;N,p):=\sum\limits_{n=1}^{N} n^{-s} + \chi(s) \sum\limits_{n=1}^{N} n^{s-1}-
\frac{(-1)^{N} }{2} \Big[ {\mathcal I}_{N+\frac{1}{2},p}(s)+ 
\chi(s) \overline{\ii}_{N+\frac{1}{2},p}(1- s)\Big], 
\end{equation}
where  $\overline{\ii}_{N+\frac{1}{2},p}(\cdot)$ is defined via \eqref{def:bar_f} and
\begin{equation}\label{def:chi}
\chi(s):=\frac{(2\pi)^s}{2\cos(\pi s/2) \Gamma(s)}.
\end{equation}
Our approximations to the Riemann zeta function are given by
\begin{equation}\label{def:zeta_p}
\zeta_p(s):=F(s;N_t,p),
\end{equation}
where $N_t=\lfloor \sqrt{t/(2\pi)} \rfloor$. 

These approximations $\zeta_p(s)$ are very easy to evaluate. They require only the precomputed values of $\{\omega_{p,j}\}_{0\le j \le p}$ and $\{\lambda_{p,j}\}_{1\le j \le p}$, along with elementary functions such as the exponential function and the logarithm function, as well as the Gamma function (which is needed to compute $\chi(s)$). The Gamma function can be efficiently computed to high precision via the Stirling series \cite{Johansson_2023}. A MATLAB implementation for computing $\zeta_8(s)$ is provided in Appendix B.
The main question now is: how closely do these functions $\zeta_p(s)$ approximate $\zeta(s)$? While we do not offer any rigorous error bounds, we aim to demonstrate the accuracy of these approximations through extensive numerical experiments.

To test the accuracy of these approximations in a vertical strip $a\le \re(s) \le b$, we aim to compute 
$$
 \max\limits_{a\le \sigma \le b} \big\vert \zeta_p(\sigma+\i t)-
\zeta(\sigma+\i t)\big\vert
$$
and plot this quantity as a function of $t$. Since computing the exact maximum is impractical, we approximate it by evaluating the function at one hundred points $\sigma_k$ equally spaced on the interval $[a,b]$. Thus, we define
$$
\Delta_p(t;a,b):=\max\limits_{0\le k \le 100} \big\vert \zeta_p(\sigma_k+\i t)-
\zeta(\sigma_k+\i t)\big\vert, \; {\textnormal{ where }} \; \sigma_k=a+(b-a)k/100.
$$
The choice of $100$ points in our definition of $\Delta_p$ is arbitrary and not essential. Using  $50$ or $200$ points would produce nearly identical numerical results. When focusing on the critical strip, we simplify the notation and write $\Delta_p(t):=\Delta_p(t;0,1)$. 

 To compute $\Delta_p(t;a,b)$, we  require benchmark values of $\zeta(s)$ computed to sufficiently high precision. These were obtained using Galway's method \cite{Galway_2001} (see Section \ref{section2} for more details). All computations were performed in Fortran 90 using David Bailey's MPFUN2020 \cite{Bailey_2020} arbitrary-precision package.

\begin{figure}[t!]
\centering
\subfloat[]{\label{fig_4a}\includegraphics[height =6.5cm]{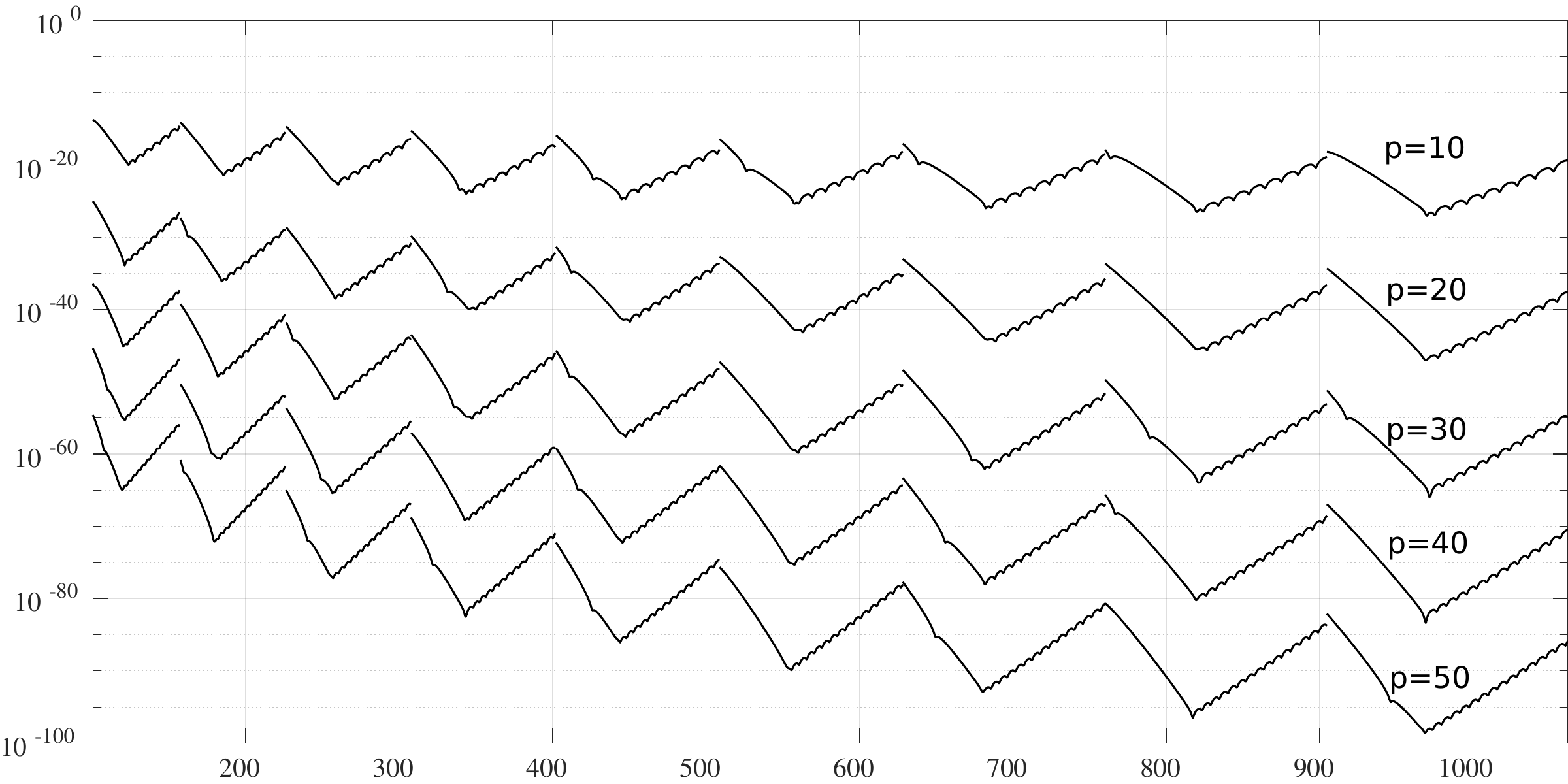}} \\
\subfloat[]{\label{fig_4b}\includegraphics[height =6.5cm]{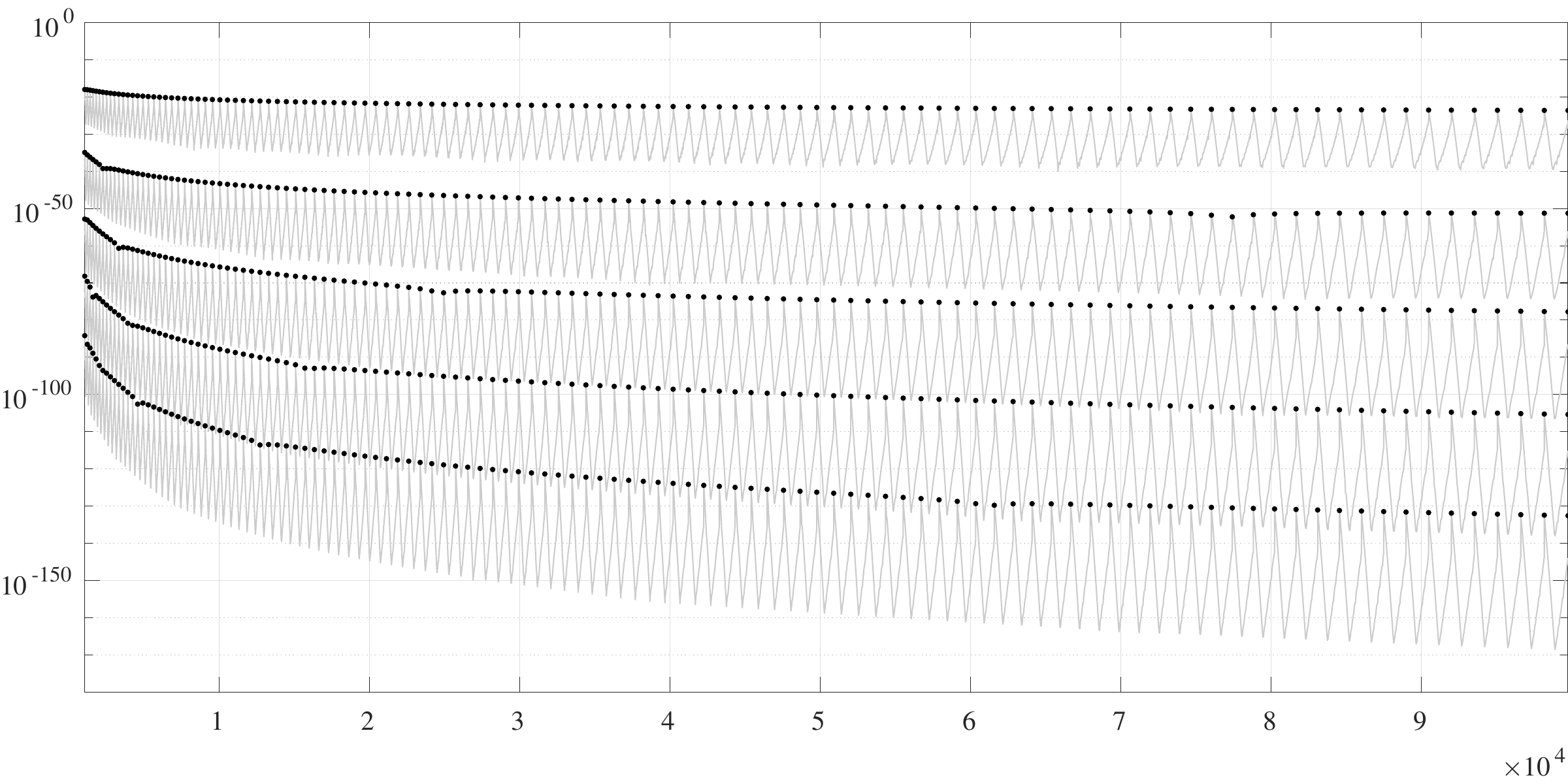}} 
\caption{The values of $\Delta_p(t):=\Delta_p(t;0,1)$ for $p \in \{10,20,30,40,50\}$. The black dots on plot (b) correspond to values of $\Delta_p(2 \pi n^2)$. } 
\label{fig4}
\end{figure}

As a first test,  we examined the accuracy of the approximations $\zeta_p(s)$ in the critical strip. We define $t_n:=2\pi n^2$, which represents the values of $t$ at which $N_t=\lfloor \sqrt{t/(2\pi)} \rfloor$  increases by one. Figure \ref{fig4} displays graphs of  $\Delta_p(t)$ over two ranges $t_4 \le t \le t_{13}$ (Figure \ref{fig_4a})
and $t_{13} \le t \le t_{126}$ (Figure \ref{fig_4b}), for $p\in \{10,20,30,40,50\}$.
We observe that $\Delta_p(t)$ exhibits a similar  pattern on each interval $[t_n, t_{n+1}]$: its values tend to be largest near the endpoints of these intervals and smallest near their midpoints.  
The numerical results presented in Figure \ref{fig4} suggest that the following error bounds hold for $s=\sigma+ \i t$ and $0\le \sigma \le 1$:
\begin{itemize}
\item $|\zeta_{10}(s)-\zeta(s)|<10^{-15}$ when $t>250$ and $|\zeta_{10}(s)-\zeta(s)|<10^{-20}$ when $t>6000$;
\item $|\zeta_{20}(s)-\zeta(s)|<10^{-30}$ when $t>350$ and
 $|\zeta_{20}(s)-\zeta(s)|<10^{-50}$ when $t>65000$;
 \item $|\zeta_{50}(s)-\zeta(s)|<10^{-100}$ when $t>4000$;
\end{itemize}

\begin{figure}[t!]
\centering
{\includegraphics[height =6.5cm]{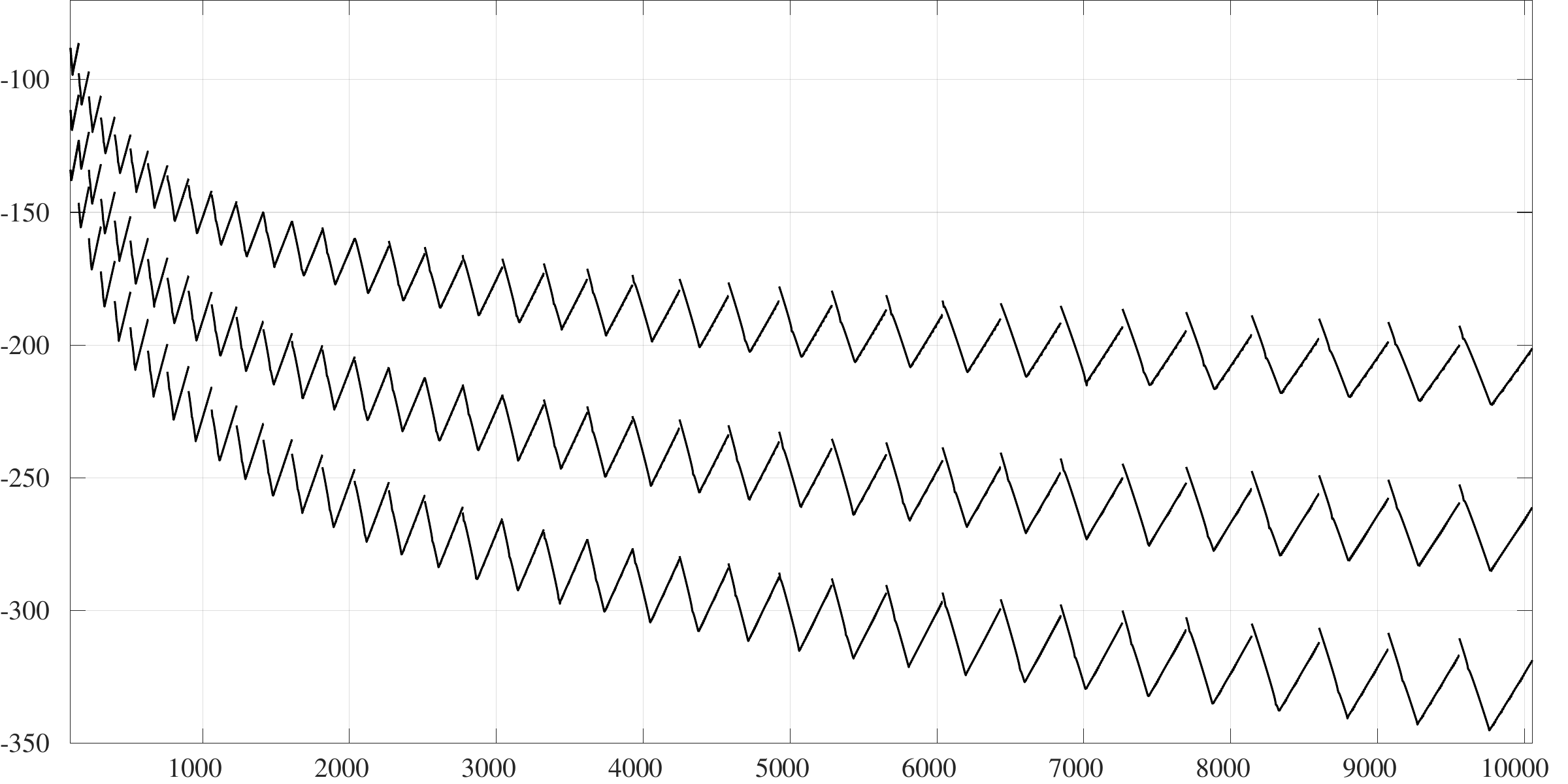}} 
\caption{The values of $\log_{10} (\Delta_p(t))$  for $p \in \{90,120,150\}$.} 
\label{fig_large_p}
\end{figure}

\begin{figure}[t!]
\centering
{\includegraphics[height =6.5cm]{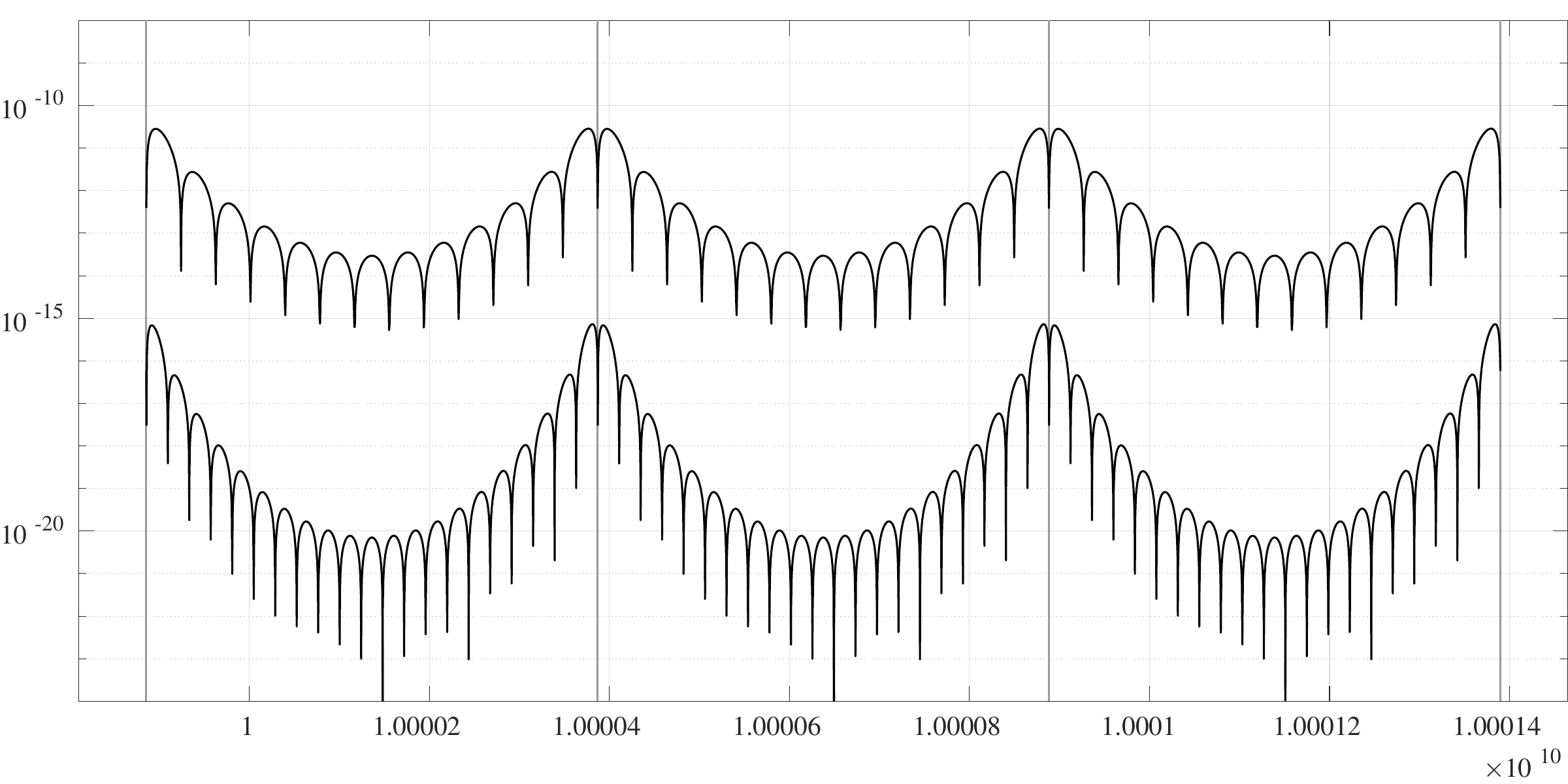}} 
\caption{The values of $|\zeta_p(1/2+\i t) - \zeta(1/2+\i t)|$ for $p \in \{3,5\}$ and $t$ close to $10^{10}$. The gray vertical lines show the locations of $t_n=2\pi n^2$ for $39894\le n \le  39897$.} 
\label{fig5}
\end{figure}

In Figure \ref{fig_large_p}, we plot $\Delta_p(t)$ for $p \in \{90, 120, 150\}$ over the range $t_4\le t\le t_{40}$. We observe that the error $|\zeta_p(s)-\zeta(s)|$ decreases significantly as $p$ increases. These numerical results suggest that for $0\le \sigma \le 1$, the following error bounds hold: 
\begin{itemize}
\item $|\zeta_{120}(s)-\zeta(s)|<10^{-200}$ when $t>1650$;
\item $|\zeta_{150}(s)-\zeta(s)|<10^{-300}$ when $t>6900$.
\end{itemize}

Next, we examine the accuracy of our approximations for very large values of $t$. In Figure \ref{fig5} we plot the values of 
$|\zeta_p(1/2+\i t) - \zeta(1/2+\i t)|$ for $p \in \{3, 5\}$ and $t$ near $10^{10}$. Again, we observe that these approximations are highly accurate: for these values of $t$, we find  $|\zeta_3(1/2+\i t) -\zeta(1/2+\i t)|\le 10^{-10}$ and the corresponding errors for $\zeta_5$ are smaller than $10^{-15}$. Interestingly, these errors exhibit certain patterns: within each interval $t\in [t_n, t_{n+1}]$, the error $|\zeta_3(1/2+\i t) -\zeta(1/2+\i t)|$ is very small at the endpoints of the interval and also at 12 equally spaced points  inside it. A similar pattern holds for $|\zeta_5(1/2+\i t) -\zeta(1/2+\i t)|$, but with 20  equally spaced points inside each interval. This is not a coincidence: in general, for very large $t$, the error 
$|\zeta_p(1/2+\i t) -\zeta(1/2+\i t)|$ is small at the endpoints of each interval $[t_n, t_{n+1}]$ and at $4p$  equally spaced points within it. From this pattern, the reader may guess  how the coefficients $\omega_{p,j}$ and $\lambda_{p,j}$ were chosen, but we defer this discussion to Section \ref{section2}.

\begin{figure}[t!]
\centering
{\includegraphics[height =6.5cm]{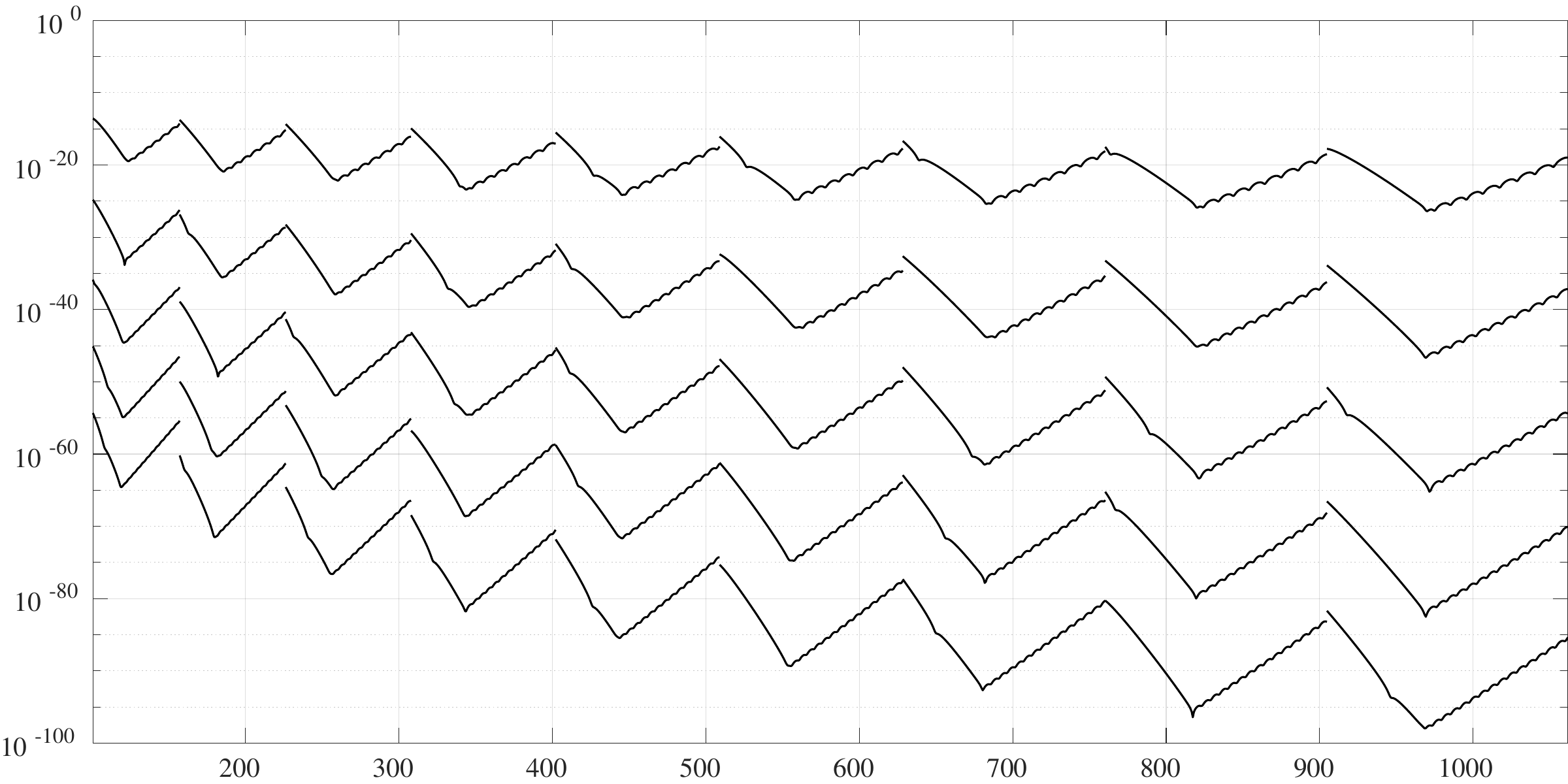}} 
\caption{The values of $\Delta^{(1)}_p(t)$  for $p \in \{10,20,30,40,50\}$. } 
\label{fig6}
\end{figure}

Computing numerical values of $\zeta'(s)$ is also important  (see, for example,  \cite{Hiary_2011b}).
Our approximations to $\zeta'(s)$ are defined as 
\begin{equation*}
\zeta^{(1)}_p(s):=
\frac{\partial}{\partial s} F(s;N,p) \Big\vert_{N=N_t},
\end{equation*}
where $F(s;N,p)$ is defined in \eqref{def:F_sNP}. 
Note that we cannot define $\zeta^{(1)}_p(s)$ as $ \frac{\d}{\d s}\zeta_p(s)$, as the latter function is not continuous (thus, not differentiable) at $t=t_n$. Figure \ref{fig6} displays graphs of
$$
\Delta^{(1)}_p(t):=\max\limits_{0\le k \le 100} \big\vert \zeta_p^{(1)}(\sigma_k+\i t)-
\zeta'(\sigma_k+\i t)\big\vert, \; {\textnormal{ where }} \; \sigma_k=k/100.
$$
The values of $\zeta'(s)$ are also computed via Galway's method \cite{Galway_2001} (details can be found in Section \ref{section2}). 
The graphs in Figure \ref{fig6} closely resemble those in Figure \ref{fig4}, suggesting that the error  $|\zeta^{(1)}_p(s)-\zeta'(s)|$ is of similar order of magnitude to $|\zeta_p(s)-\zeta(s)|$.

\begin{figure}[t!]
\centering
\subfloat[]{\label{fig_7a}
\includegraphics[height =6.35cm]{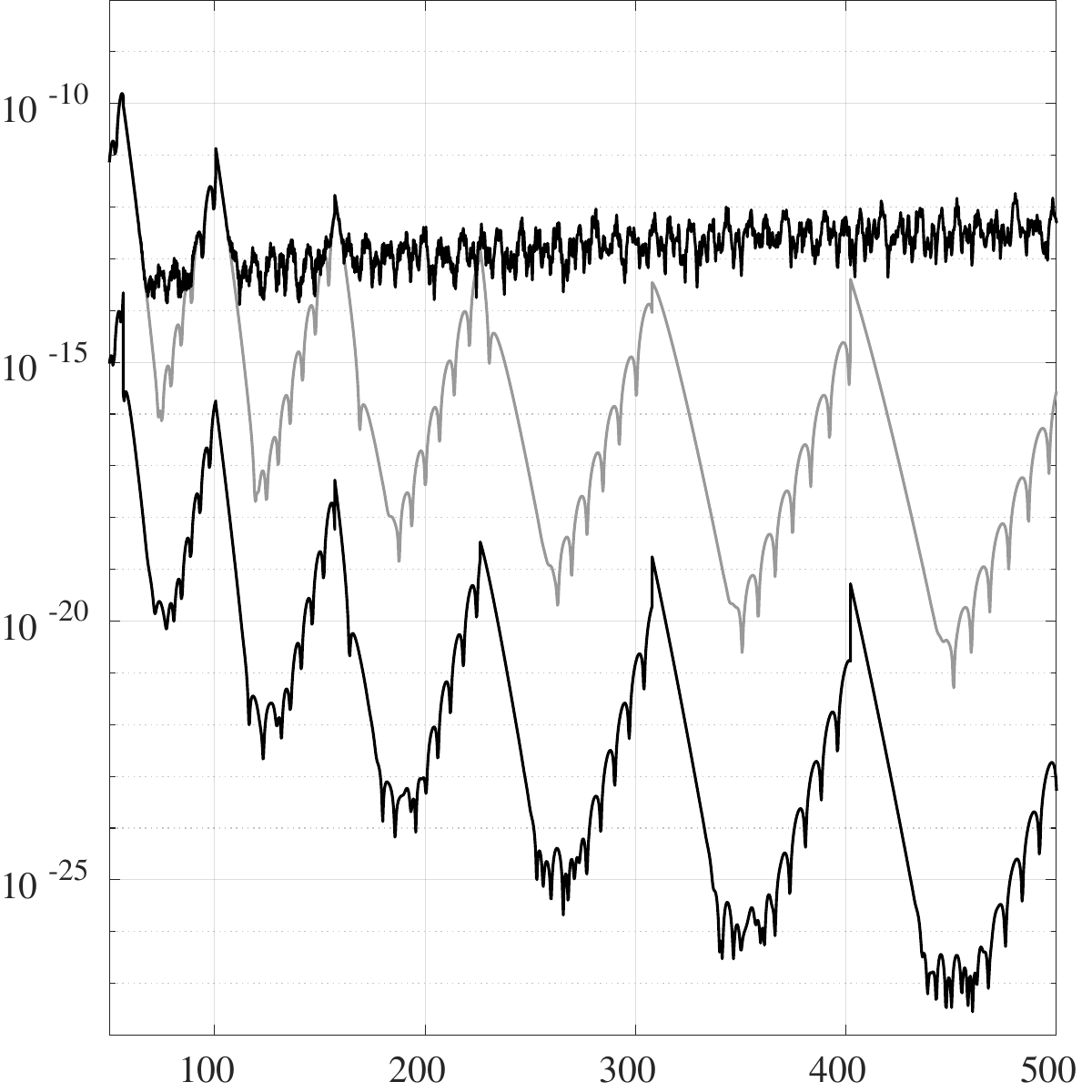}}
\subfloat[]{\label{fig_7b}
\hspace{0.2cm}
\includegraphics[height =6.5cm]{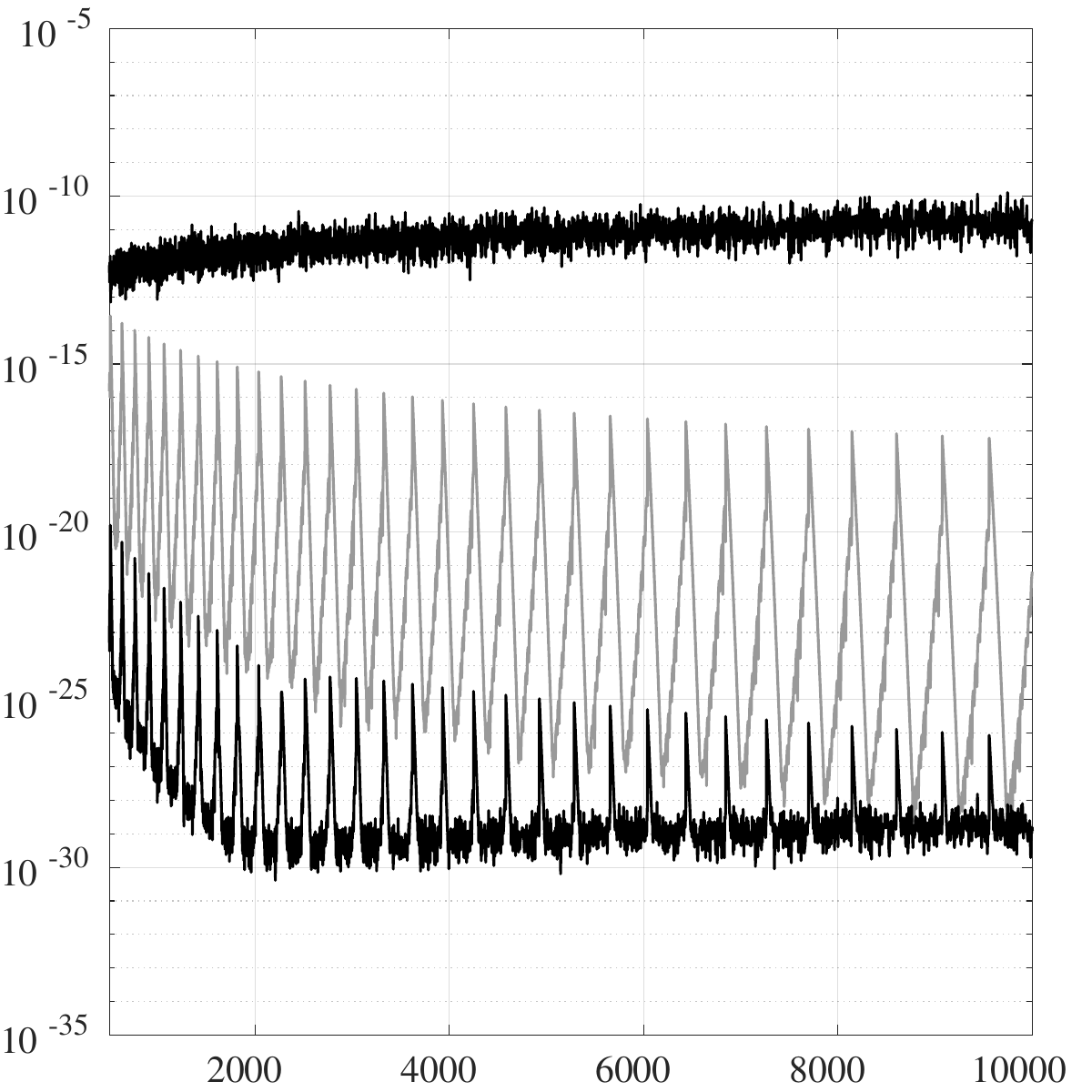}}
\caption{The values of $\Delta_p(t;1/2,2)$. The top curve corresponds to $p=8$ implemented in double precision.
The middle gray curve corresponds to $p=8$
 and the bottom curve to $p=12$, both implemented in quadruple precision.} 
\label{fig7}
\end{figure}

The computations in the previous examples were performed with sufficiently high precision to ensure that we could analyze the approximation error $\zeta_p(s)-\zeta(s)$ without having to worry about rounding errors. A natural question then arises: how accurate are these approximations when implemented in double (or quadruple) precision? To investigate this, we implemented $\zeta_8(s)$ and $\zeta_{12}(s)$ in double and quadruple precision in Fortran 90 and MATLAB. The corresponding code can be downloaded at \href{https://kuznetsovmath.ca/}{kuznetsovmath.ca/}. We tested the accuracy of these approximations in the strip $1/2\le \re(s) \le 2$. This particular strip was chosen because $\zeta(s)$ does not grow too rapidly as $t \to +\infty$ there, unlike in any strip where $\re(s)<1/2$. The results of these computations are presented in Figure \ref{fig7}. 

The top graph shows the error $\Delta_8(t;1/2,2)$, where $\zeta_8(s)$ was implemented in double precision (the benchmark values $\zeta(s)$ were computed in higher precision). We observe that rounding errors dominate the approximation error for $t>200$. These rounding errors primarily arise when computing $\chi(s)$ and evaluating the values of $n^{-s}$ and $n^{s-1}$ for $s=\sigma+\i t$ with even moderately large $t$. Indeed, it is easy to verify numerically that when computing the value of $2^{\i t}$ in double precision for $t$ of the order $10^2$, the precision of the result is around $10^{-14}$, and when $t$ increases to $10^3$, the precision drops to about $10^{-13}$. This suggests that roughly one decimal digit of precision is lost each time 
$t$ increases by a factor of ten. Precision is also lost due to cancellation errors when adding many terms in the main sum in \eqref{def:F_sNP}, though this effect likely plays a lesser role compared to rounding errors.   

The middle graph in Figure \ref{fig7} (shown in gray) displays  the errors $\Delta_8(t;1/2,2)$, where this time  $\zeta_8(s)$ was implemented in quadruple precision. We observe that rounding errors do not affect the results over this range of $t$. The quadruple precision implementation of $\zeta_8(s)$ produces errors smaller than $10^{-13}$ for $t>250$ and smaller than $10^{-15}$ for $t>2000$ in the strip $1/2 \le \sigma \le 2$. Of course, when $t$ becomes very large (of the order of $10^{10}$ or greater), rounding errors will eventually become noticeable. 

The bottom graph in Figure \ref{fig7} shows the errors 
$\Delta_{12}(t;1/2,2)$ for the quadruple precision implementation of $\zeta_{12}(s)$. These errors are significantly smaller: we find that  $\Delta_{12}(t;1/2,2)< 10^{-25}$ for $t>5000$. The effects of rounding errors become clearly noticeable for $t>2000$, though they do not pose a significant issue in this range, as the maximum approximation error remains larger than the rounding errors.  However, when $t$ reaches $10^{5}$ or greater, rounding errors will dominate the approximation error, and the graph will resemble the top graph (though with a much smaller overall error, of order $10^{-27}$).

This paper is organized as follows. In Section \ref{section2} we  discuss how our approximations $\zeta_p$ were derived, and we present an algorithm for computing the  coefficients $\omega_{p,j}$ and $\lambda_{p,j}$. In Section \ref{Section_compare_with_RS} we compare our approximations with the classical Riemann-Siegel approximations, and in Section \ref{Section_computing_zeta} we discuss how the approximations $\zeta_p$ are used in the algorithm for computing $\zeta(s)$ for arbitrary complex $s$. 

\section{Deriving the approximations and computing $\omega_{p,j}$ and $\lambda_{p,j}$}\label{section2}

We begin with the following result: for every integer $N\ge 0$ and $s\in \c$
\begin{equation}\label{eqn:zeta_as_RS}
\zeta(s)=\sum\limits_{n=1}^N n^{-s} + \chi(s) \sum\limits_{n=1}^N n^{s-1}-
\frac{(-1)^{N} }{2} \Big[ {\mathcal I}_{N+\frac{1}{2}}(s)+ 
\chi(s) \overline{\ii}_{N+\frac{1}{2}}(1- s)\Big],
\end{equation}
where 
\begin{equation}\label{def:I_M}
{\mathcal I}_M(s):= \int_{\r} \frac{e^{-\pi x^2-2\pi M  \theta  x}}{\cosh(\pi \theta x)} 
 \Big(M+\frac{x}{\theta}\Big)^{-s} \d x, 
 \end{equation} 
 and $\theta:=\exp(-\pi \i /4)$. 
The function $\overline{\ii}_M$ is defined via \eqref{def:bar_f}.  
The above result is equivalent to formulas (1.1), (1.3) and (1.4) in \cite{Galway_2001}, after a change of variable of integration 
$z=N+1/2+\theta x$ in \cite{Galway_2001}[equation (1.1)]. It can also be derived from formulas  (1.1) and (3.2) in
\cite{Reyna_2011}.

Formula \eqref{eqn:zeta_as_RS} was used by Galway \cite{Galway_2001} to derive a numerical quadrature algorithm for computing $\zeta(s)$. We used a special case of Galway's approximations to compute the benchmark values of $\zeta(s)$ and $\zeta'(s)$. We approximate $\ii_M(s)$ by 
\begin{equation}\label{def_ii_M_sh}
\ii_M(s;h):=h \sum\limits_{k \in {\mathbb Z}} 
\frac{e^{-\pi (kh)^2-2\pi M  \theta  kh}}{\cosh(\pi \theta kh)} 
 \Big(M+\frac{kh}{\theta}\Big)^{-s},
\end{equation} 
where $h>0$. Note that for $M\ge 1/2$, the integrand in \eqref{def:I_M} is analytic in a strip $|\im(x)|<\pi/\sqrt{2}$ and decays very rapidly as $|x| \to \infty$. 
According to Lemma 2.1 and estimates (2.2) and (2.3) in \cite{Galway_2001}, for any $s\in \c$, $M\ge 1/2$ and $c\in (0,\pi/\sqrt{2})$  we have
\begin{equation}\label{I_M_sh_error}
\ii_M(s;h)-\ii_M(s)=O\big(\exp(-c /h)\big), \;\;\; h\to 0^+. 
\end{equation}
This result also follows from \cite{Trefethen2014}[Theorem 5.1]. The implied constant in the  big-O term depends on $s$, $M$ and $c$.

\begin{figure}[t]
\centering
{\includegraphics[height =6.5cm]{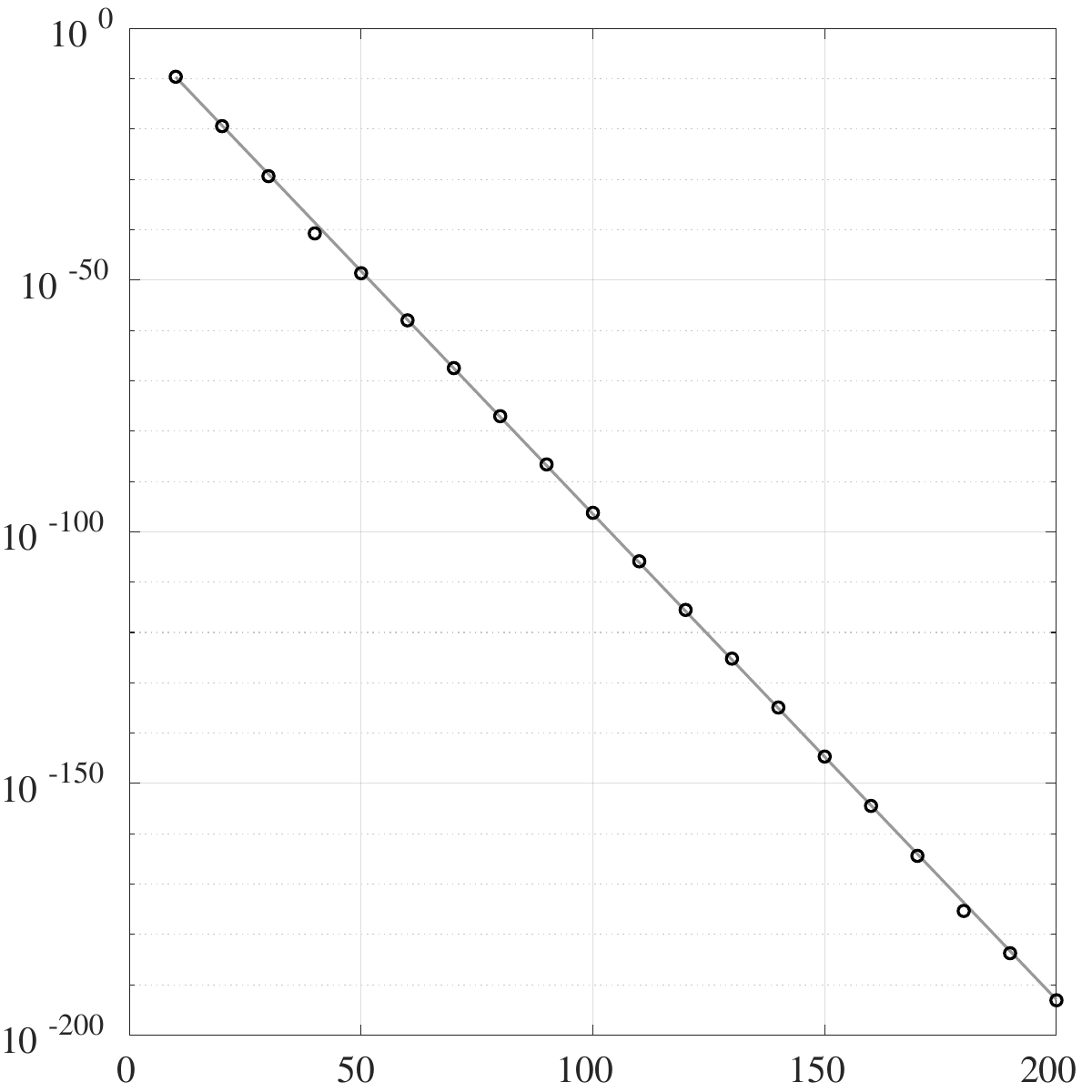}} 
\caption{The values of $|\zeta(1/2+\i \gamma_{30};h)|$ (circles) and $\exp(-\frac{\pi}{\sqrt{2}} \times \frac{1}{h})$ (gray line), where $\frac{1}{h}$ is on the $x$-axis.  }
\label{fig8}
\end{figure}

We define for $N\ge 0$ 
$$
G(s;N,h):=\sum\limits_{n=1}^N n^{-s} + \chi(s) \sum\limits_{n=1}^N n^{s-1}-
\frac{(-1)^{N} }{2} \Big[ {\mathcal I}_{N+\frac{1}{2}}(s;h)+ 
\chi(s) \overline{\ii}_{N+\frac{1}{2}}(1- s;h)\Big].
$$
Formula \eqref{I_M_sh_error} implies that for any non-negative integer $N$ we have $G(s;N,h) \to \zeta(s)$ as $h \to 0^+$, so that we have an infinite family of approximations to $\zeta(s)$, allowing us to choose (for every $s$) the one that is easiest to compute. It is well known \cite{Galway_2001,Titchmarsh1987} that the optimal choice is  $
N=N_t=\lfloor \sqrt{t/(2\pi)} \rfloor$ (recall that $s=\sigma+\i t$ and $t>0$). This choice ensures that the saddle point of the integrand in \eqref{def:I_M} is as close as possible to the real line, so that the resulting integral is easy to evaluate numerically.
We define $\zeta(s;h)=G(s;N_t,h)$. This is  how we compute the benchmark values of $\zeta(s)$ to high precision. Galway \cite{Galway_2001}  developed  more general (and more efficient) approximations, but this simpler version suffices for our purposes. To implement $\zeta(s;h)$, we used Stirling's series algorithm for computing the Gamma function \cite{Johansson_2023} and truncated the infinite sum in \eqref{def_ii_M_sh} once the terms became smaller than $10^{-d}$, where $d$ is the working precision. To illustrate the accuracy of this approximation, we conducted the following experiment. We downloaded the value of the 30th zero $1/2+\i \gamma_{30}$ of the Riemann zeta function on the critical line, which was computed to 1000 decimal digits in \cite{Odlyzko_tables}, and we evaluated $\zeta(1/2+\i \gamma_{30};h)$ for a decreasing sequence of $h$. The results, shown in Figure \ref{fig8}, demonstrate that the  values of $\zeta(1/2+\i \gamma_{30};h)$ decrease precisely at the rate $\exp(-\frac{\pi}{\sqrt{2}} \times \frac{1}{h})$, which is consistent with \eqref{I_M_sh_error}. 

To compute the benchmark values of $\zeta'(s)$ (needed for the results in Figure \ref{fig6}), we approximated it using
$$
\zeta^{(1)}(s;h):=\frac{\partial}{\partial s} G(s;N;h) \Big \vert_{N=N_t}. 
$$
We also conducted several numerical experiments confirming that $\zeta^{(1)}(s;h) \to \zeta'(s)$ as $h\to 0^+$ at the same rate $\exp(-\frac{\pi}{\sqrt{2}} \times \frac{1}{h})$. 

We now explain the ideas behind our approximations and the method used to compute the coefficients $\omega_{p,j}$ and $\lambda_{p,j}$. Formula \eqref{eqn:zeta_as_RS} shows that we can compute $\zeta(s)$ once we have an effective way to evaluate  $\ii_M(s)$.
We rewrite \eqref{def:I_M} in the form
\begin{equation}\label{def:I_M2}
{\mathcal I}_M(s):= M^{-s} \int_{\r}  e^{g(x;s)} \times \frac{e^{-\pi x^2}}{\cosh(\pi \theta x)}  
 \d x=M^{-s} \int_{\r} e^{g(x;s)} \eta(\d x),
 \end{equation}  
where we denoted
$$
g(x)=g(x;s):=-2 M \pi \theta x - s \ln\Big( 1  + \frac{x}{M\theta} \Big) \;\;\;\; {\textnormal{ and }} \;\;\;\; 
\eta(\d x):=\frac{e^{-\pi x^2}}{\cosh(\pi \theta x)}  \d x.
$$
We consider the complex numbers $\omega_{p,j}$ and $x_{p,j}=\lambda_{p,j}/\theta$  as the weights and nodes of a quadrature
with respect to the complex-valued measure $\eta(\d x)$
\begin{equation}\label{eqn:quadrature_condition}
\int_{\r} f(x)  \eta(\d x) \approx 
\omega_{p,0} f(0)+  \sum\limits_{j=1}^p  \omega_{p,j} \Big[f(x_{p,j}) +
f(-x_{p,j})\Big].  
\end{equation}
A classical way to determine  the coefficients $\omega_{p,j}$ and $x_{p,j}$ is to require that \eqref{eqn:quadrature_condition} holds exactly for a chosen set of functions $f_k(x)$. 
For example, the Gaussian quadrature method requires \eqref{eqn:quadrature_condition} to be exact for all polynomials of degree $\le 4p+1$. In fact, this was the first approach we explored, and while it produces decent results, the method we introduce next provides better accuracy and is somewhat easier to implement.

For the remainder of this section, we assume that $s$ lies within a fixed vertical strip $a\le \sigma \le b$ and we set $M=N_t+1/2$, where $s=\sigma+\i t$, $t>0$ and $N_t=\lfloor \sqrt{t/(2\pi)} \rfloor$. We now examine the asymptotic behavior of the integrand in \eqref{def:I_M2} as $t\to +\infty$. 
Expanding $g(x;s)$ in a Taylor series around $x=0$, we obtain 
for $x\in \c$ with $|x|<M$
\begin{equation}\label{eqn:g_Taylor}
g(x;s)=\Big(-\frac{\i s}{2\pi M}-M\Big) 2\pi \theta x+  
\frac{1}{2} \frac{s}{M^2 \theta^2 } x^2 + \sum\limits_{k\ge 3} 
\frac{(-1)^k}{k} \frac{s }{M^k \theta^k} x^k.
\end{equation}
The first coefficient in this expansion simplifies to
$$
-\frac{\i s}{2\pi M}-M=\frac{t-\i \sigma}{2 \pi (N_t + \frac{1}{2})} - N_t  - \frac{1}{2}=\frac{t}{2\pi N_t}-N_t-1+O(t^{-1/2}),  
$$
as $t\to +\infty$. 
Defining 
\begin{equation*}
B(t):=\frac{t}{2\pi N_t}-N_t-1. 
\end{equation*}
and using the inequality $N_t \le \sqrt{t/(2\pi)} < N_t+1$, we obtain the bound
\begin{equation}\label{bounds_Bt}
-1\le B(t)\le 1+O(t^{-1/2}).
\end{equation}
The second coefficient in \eqref{eqn:g_Taylor} simplifies to
$$
\frac{1}{2} \frac{ s}{M^2 \theta^2 } = - \frac{t}{2 (N_t+\frac{1}{2})^2}+O(t^{-1})=-\pi+O(t^{-1/2}). 
$$
When $t$ is large enough, we have $|s|<2t$ and $M=N_t+1/2 > \sqrt{t}/3$ (since $\sqrt{2\pi}<3$). Therefore, for $t$ large enough, 
the sum in \eqref{eqn:g_Taylor} can be estimated as follows
\begin{equation}
\Big| \sum\limits_{k\ge 3} 
\frac{(-1)^k}{k} \frac{s }{M^k \theta^k} x^k \Big| < 
2t  \sum\limits_{k \ge 3} \Big(3 |x|/\sqrt{t} \Big)^{k}
=  \frac{ 54 |x|^3 }{\sqrt{t}- 3|x| },
\end{equation}
for $x\in \c$ with $|x|<\sqrt{t}/3$. 
Thus, we conclude that  
\begin{equation}\label{eqn:g_asymptotics}
g(x;s)=- \pi x^2+2\pi \theta B(t) x  + O(t^{-1/2}),
\end{equation}
as $t \to +\infty$, uniformly in $x$ on compact subsets of $\c$. 

\begin{figure}[t]
\centering
{\includegraphics[height =6.5cm]{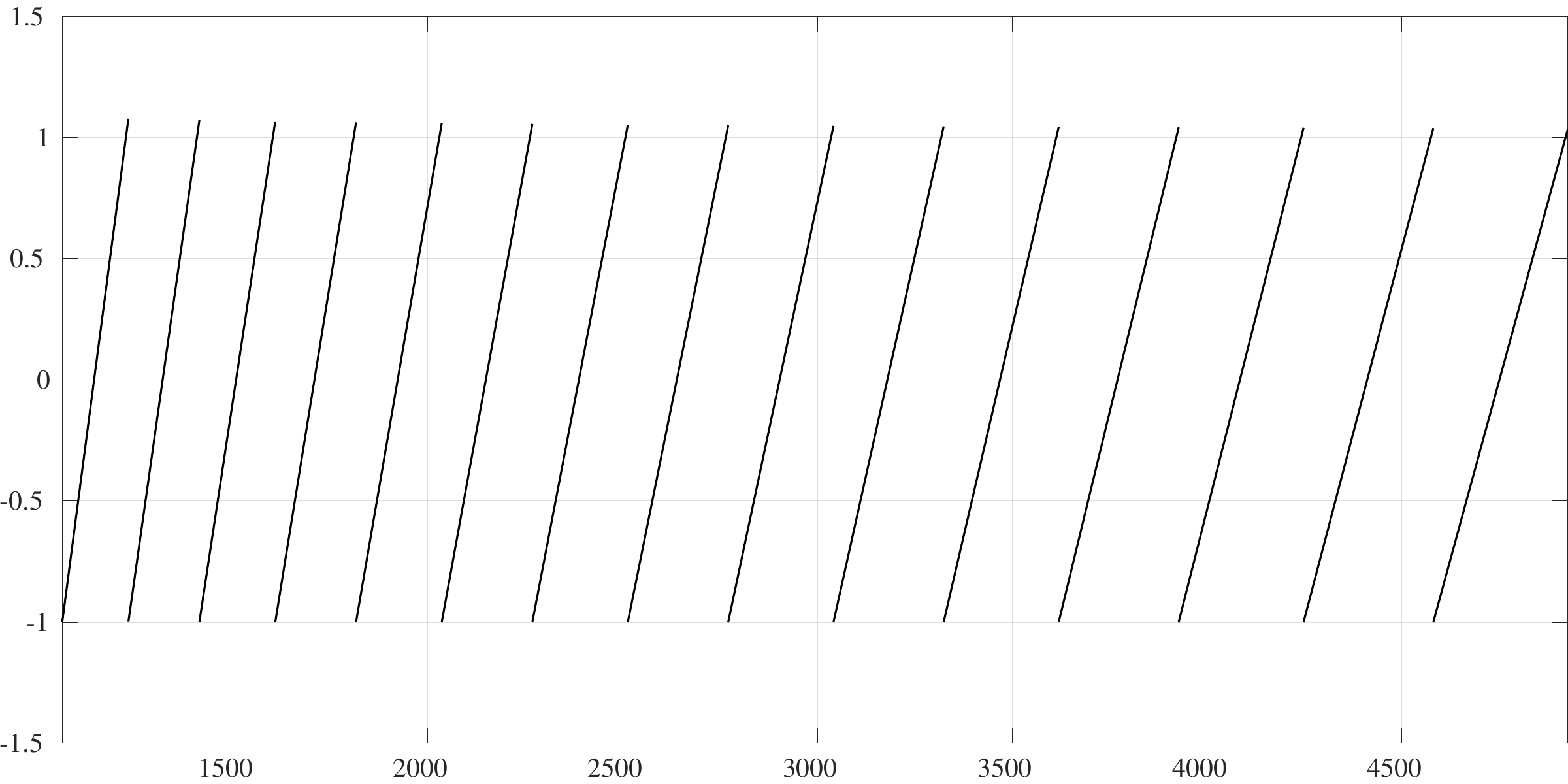}} 
\caption{The graph of  $B(t)$.  }
\label{fig9}
\end{figure}

For $t$ large, the integral in \eqref{def:I_M2} is therefore approximately equal to 
$$
\int_{\r} 
e^{-\pi x^2+2\pi \theta B(t) x} \eta(\d x), 
$$
where $B(t)$ typically lies in the interval $[-1,1]$, though occasionally it may slightly exceed one (see Figure \ref{fig9}). 
To determine the quadrature coefficients  $\{\omega_{p,j}\}_{0\le j \le p}$ and 
$\{x_{p,j}\}_{1\le j \le p}$, we require that   \eqref{eqn:quadrature_condition} holds exactly for a set of $4p+2$ functions 
$$
f_k(x)=e^{-\pi x^2 + 2\pi \theta y_{p,k} x}, \;\;\;  k=0,1,\dots,4p+1,
$$
where the points  $\{y_{p,k}\}_{0\le k \le 4p+1}$ are equally spaced in the interval $[-1,1]$:
$$
y_{p,k}:=-1+\frac{2k}{4p+1}.
$$
Thus, we arrive at the following problem: we seek $2p+1$ complex numbers $\{\omega_{p,j}\}_{0\le j \le p}$ and 
$\{x_{p,j}\}_{1\le j \le p}$ such that 
\begin{equation}\label{eqn:quadrature_condition2}
\int_{\r}
 \frac{e^{-2\pi x^2+2\pi \theta y_{p,k} x}}{\cosh(\pi \theta x)}  \d x  =
\omega_{p,0} +  \sum\limits_{j=1}^p  \omega_{p,j} e^{-\pi x_{p,j}^2}  \Big[e^{2\pi \theta y_{p,k} x_{p,j}} +e^{-2\pi \theta y_{p,k} x_{p,j}}\Big],
\end{equation}
for all  $k=0,1,\dots,4p+1$. 

To simplify the above problem, we note that the integral on the left-hand side of \eqref{eqn:quadrature_condition2} is one of the Mordell integrals \cite{Mordell2} and can be computed explicitly. Using formula (7) in \cite{Kuznetsov2007} (or transforming the integral in \cite{Gabcke1979}[Theorem 4.1.1]) 
we evaluate 
\begin{equation}\label{eqn:Mordell_integral}
H(y):=\int\limits_{\r} \frac{e^{-2\pi x^2+ 2 \pi  \theta  x y }}{\cosh(\pi \theta x)}   \d x=
\frac{1}{\cos(\pi y)} \Big[\sqrt{2} \cos(\pi y/2) e^{-\frac{\pi \i}{8}
(4y^2+1)} -e^{-\frac{\pi \i}{4}}\Big].
\end{equation}
Note that $H(y)$ is an entire function of $y$. 
We denote $\mu_{p,k}:=H(y_{p,k})$ 
and introduce new variables $\{u_j,z_j\}_{-p \le j \le p}$, defined as follows: $u_0=\omega_{p,0}$, $z_0=1$ and for $j=1,2,\dots,p$
\begin{equation*}
u_j:=\omega_{p,j} e^{-\pi x_{p,j}^2 - 2\pi \theta x_{p,j}}, \;\;\;
u_{-j}:=\omega_{p,j} e^{-\pi x_{p,j}^2 + 2\pi \theta x_{p,j}}, \;\;\;
z_j:=e^{2\pi \theta \frac{2}{4p+1} x_{p,j}}, \;\;\; 
z_{-j}:=1/z_j. 
\end{equation*}
Note that $u_j$ and $z_j$ also depend on $p$, but we suppress the dependence on $p$ to simplify  notation. 
In these new variables, the system of equations \eqref{eqn:quadrature_condition2} simplifies to 
\begin{equation}\label{eqn:Gaussian_quadrature}
\sum\limits_{j=-p}^p u_j z_j^k = \mu_{p,k}, \;\;\; k=0,1,\dots,4p+1. 
\end{equation}
Our task is now to determine the values of $\{u_j,z_j\}_{-p \le j \le p}$ satisfying $z_0=1$, $z_{-j}=1/z_j$ and $u_{-j}=z_j^{4p+1} u_j$ for $j=1,2,\dots,p$, which solve the system of equations \eqref{eqn:Gaussian_quadrature}.

The solution to the above system of equations is obtained using classical methods of Gaussian quadrature \cite{GONZALEZV1995,Laurie_2001,MILOVANOVIC2005,Pozza2016}. Let ${\mathscr P}_n$ be the space of polynomials with complex coefficients of degree not exceeding $n$. Define a linear functional $L$ acting on ${\mathscr P}_{4p+1}$ via
$$
L[x^k]=\mu_{p,k}, \;\;\; k=0,1,\dots,4p+1. 
$$
The system \eqref{eqn:Gaussian_quadrature} is equivalent to the condition
\begin{equation}\label{eqn:Gaussian_quadrature3}
L[Q]=\sum\limits_{j=-p}^p u_j Q(z_j) \;\;
{\textnormal{ for all }} \; Q \in {\mathscr P}_{4p+1}. 
\end{equation}
Thus, we can recognize $z_j$ and $u_j$ as the nodes and weights of the Gaussian quadrature for the linear functional $L$. 
We find these nodes and weights via orthogonal polynomials with respect to the linear functional $L$ (see \cite[Definition 1.1]{Pozza2016}). We denote $m:=2p+1$, and we aim to find polynomials $\{P_{m,n}\}_{0\le n \le m}$ satisfying ${\textnormal{deg}} ( P_{m,n})=n$ and the orthogonality condition: for $1 \le n \le m$  
$$
L[P_{m,n} Q]=0 \; {\textnormal{ for all }} \; Q\in {\mathscr P}_{n-1}. 
$$
We will suppress the dependence on $m$ in what follows and will write simply $P_n=P_{m,n}$.

It is known (see \cite[Theorem 3.2]{Pozza2016}) that these orthogonal polynomials exist if and only if $L$ is {\it quasi-definite} on ${{\mathscr P}}_m$. The latter condition means that certain Hankel determinants involving the coefficients $\mu_{p,k}$ are non-zero (see formula (2.1) and Definition 3.1 in \cite{Pozza2016}). Assuming that the orthogonal polynomials $P_n$ exist, they can be computed   via the three-term recurrence relation 
\begin{equation}\label{three_term_recurrence}
 P_{n+1}(x)=(x-a_n) P_n(x) - b_n P_{n-1}(x), \;\;\; n=0,1,\dots m-1, 
\end{equation}
 with initial conditions $P_{-1}(x)\equiv 0$, $P_0(x) \equiv 1$. The  recurrence  coefficients  are given by
 $$
 a_n = \frac{L[ x P_n^2]}{L[P_n^2]}, \;\;\; b_n=\frac{L[ P_n^2]}{L[P_{n-1}^2]}. 
 $$
 The above formulas can be found in equations (3.1) and (3.2) in \cite{Pozza2016} or in \cite[Theorem 2.9]{MARCELLAN2001}.

Before we proceed, we establish the following result:   
\begin{proposition}\label{prop:P_m} Assume that there exist polynomials $\{P_n\}_{0\le n \le m}$ that are orthogonal with respect to the linear functional $L$. Then 
the polynomial $P_m$  must satisfy
$P_m(x)=-x^m P_m(1/x)$ for all $x\neq 0$.
\end{proposition} 
\begin{proof} We begin by stating a preliminary result about determinants. Let $B=\{b_{i,j}\}_{1 \le i,j \le n}$ be an $n\times n$ matrix. Following \cite{Golyshev2007}, we denote by $B^{\tau}$ the transpose of  $B$  with respect to the anti-diagonal. In other words, $B^{\tau}=\{\tilde b_{i,j}\}_{1 \le i,j \le n}$ where $\tilde b_{i,j}=b_{n+1-j,n+1-i}$ for $1\le i,j \le n$. 
It is true that
\begin{equation}\label{eqn:det_B_tau}
{\textnormal{det}}(B) = {\textnormal{det}}(B^{\tau}).
\end{equation} 
This follows from the factorization $B^{\tau}={\textrm{J}} B^T {\textrm{J}}$ (see \cite{Golyshev2007}), where ${\textrm{J}}$ is an $n\times n$ exchange matrix, defined via ${\textrm{J}}_{i,j}=1$ if $i+j=n+1$ and ${\textrm{J}}_{i,j}=0$ otherwise. 

Now we can prove Proposition \ref{prop:P_m}. Since $H(y)$ is an even function (see \eqref{eqn:Mordell_integral}),   the moments  satisfy $\mu_{p,k}=\mu_{p,4p+1-k}$ for $0\le k \le 4p+1$. Using this fact and the well-known representation of orthogonal polynomials as determinants (see formula (2.2) in \cite{MARCELLAN2001}), we obtain
\begin{equation}\label{eqn:P_m_determinant}
P_m(x)= C \times  {\textnormal{det}}\begin{bmatrix}
\mu_0 & \mu_1 & \mu_2 & \cdots & \mu_{2p-1} &  \mu_{2p} & \mu_{2p} \\
\mu_1 & \mu_2 & \mu_3 & \cdots & \mu_{2p} & \mu_{2p} & \mu_{2p-1} \\
\mu_2 & \mu_3 & \mu_4 & \cdots & \mu_{2p} & \mu_{2p-1} & \mu_{2p-2} \\
\cdots  & \cdots  & \cdots  &   \cdots & \cdots & \cdots  & \cdots  \\
\mu_{2p-1} & \mu_{2p} & \mu_{2p}  & \cdots &\mu_{3} &\mu_{2} & \mu_{1} \\
\mu_{2p} & \mu_{2p} & \mu_{2p-1}  & \cdots &\mu_{2} &\mu_{1} & \mu_{0} \\
1 & x & x^2 & \cdots & x^{2p-1} & x^{2p} & x^{2p+1} 
\end{bmatrix}
\end{equation}
for some constant $C$. In the above formula, we used the simplified notation $\mu_k=\mu_{p,k}$. 
For $k=0,1,\dots,m$ we denote by $A_k$ the submatrix obtained from the matrix in \eqref{eqn:P_m_determinant} by removing the last row and the $(k+1)$-st column.  Performing Laplace expansion along the last row, we obtain
$$
P_m(x)=C \sum\limits_{k=0}^{m} (-1)^k {\textnormal{det}}(A_k) x^k. 
$$
It is straightforward  to verify that $A_k=A_{m-k}^{\tau}$. The desired result follows from \eqref{eqn:det_B_tau}.   
\end{proof}

Next, we assume that we have been able to compute orthogonal polynomials 
$\{P_n\}_{0\le n \le m}$ via the three-term recurrence relation 
\eqref{three_term_recurrence}, and we assume further that all roots $\{z_j\}_{-p \le j \le p}$ of the polynomial $P_m$ are simple. While this would be guaranteed if the moments $\mu_{p,k}$ were moments of a positive measure (which is not the case here), in all our computations we found that the roots of $P_m$ were indeed simple. Proposition \ref{prop:P_m} implies that one of the roots must be equal to $1$ and the other roots come in pairs $z_j$ and $1/z_j$. This allows us to order the roots $\{z_j\}_{-p \le j \le p}$ in increasing order of their absolute values (that is, $|z_j|\le |z_{j+1}|$) and in such a way that $z_{-j}=1/z_j$ and $z_0=1$. 
The weights of the Gaussian quadrature are computed using \cite[formula (3)]{Laurie_2001}:
 \begin{equation}
 u_j= \frac{L[P_{m-1}^2]}{ P_{m-1}(z_j) P_{m}'(z_j)}, \;\;\; -p \le j \le p.
 \end{equation}
Note that $P_{m}'(z_j)\neq 0$ due to our assumption that the roots are simple. Additionally, $P_{m-1}(z_j)\neq 0$, since otherwise, the three-term recurrence \eqref{three_term_recurrence} would imply that $P_n(z_j)=0$ for all $0\le n \le m$, which is impossible (recall that $P_0(z)\equiv 1$). 

We claim that the weights and nodes of this Gaussian quadrature  satisfy the identity $u_{-j}=z_j^{4p+1} u_j$. This follows from \eqref{eqn:Gaussian_quadrature} and the symmetry relation $\mu_{p,k}=\mu_{p,4p+1-k}$. Indeed, after we have found $z_j$ and established that $z_{-j}=1/z_j$, we can interpret \eqref{eqn:Gaussian_quadrature} as a system of linear equations in the $2p+1$ unknowns $\{u_j\}_{-p \le j \le p}$. This system has a unique solution, since its coefficient matrix is the Vandermonde matrix constructed from distinct numbers $\{z_j\}_{-p \le j \le p}$. By changing indices $k \to 4p+1-k$ and $j \mapsto -j$ in \eqref{eqn:Gaussian_quadrature}, and using the identities $\mu_{p,k}=\mu_{p,4p+1-k}$ and $1/z_{j}=z_{-j}$, we conclude that the numbers  
$\{u_{-j}z_j^{-4p-1}\}_{-p \le j \le p}$ also solve the same system of equations. By the uniqueness of the solution, it follows that $u_{-j}=u_j z_{j}^{4p+1}$ for all $-p \le j \le p$.

After we have found the nodes  $\{z_j\}_{-p\le j \le p}$ and the weights $\{u_j\}_{-p\le j \le p}$ of Gaussian quadrature \eqref{eqn:Gaussian_quadrature3}, we compute $\omega_{p,j}$ and $x_{p,j}$ from \eqref{def:tilde_w_z}. We set $\omega_{p,0}=u_0$ and
\begin{equation}\label{def:tilde_w_z}
x_{p,j}=\frac{4p+1}{4\pi \theta} \log(z_j), \;\;\;
\omega_{p,j}=u_j e^{\pi x_{p,j}^2 + 2\pi \theta x_{p,j}}, \;\;\;
 j=1,2,\dots,p.
\end{equation}
By construction, these numbers $\{\omega_{p,j}\}_{0 \le j \le p}$ and 
$\{x_{p,j}\}_{1 \le j \le p}$ satisfy the system of equations \eqref{eqn:quadrature_condition2}. 

With the coefficients $\{\omega_{p,j}\}_{0\le j \le p}$ and 
$\{x_{p,j}\}_{1\le j \le p}$ computed, we can now use \eqref{eqn:quadrature_condition} and approximate 
\begin{align*}
{\mathcal I}_M(s)&= \int_{\r} \frac{e^{-\pi x^2-2\pi M  \theta  x}}{\cosh(\pi \theta x)} 
 \Big(M+\frac{x}{\theta}\Big)^{-s} \d x\\
 &\approx \omega_{p,0} M^{-s}   + \sum\limits_{j=1}^p 
\omega_{p,j} \Big[
e^{-2 \pi M  \theta  x_{p,j}} 
 \Big(M+  \frac{x_{p,j}}{\theta}\Big)^{-s}
+e^{2  \pi M   \theta  x_{p,j}} 
 \Big(M - \frac{x_{p,j}}{\theta}\Big)^{-s} \Big].
\end{align*}
We recognize that the right-hand side is precisely $\ii_{M,p}(s)$ defined in \eqref{def:ii_Mp} (recall that we denoted $\theta x_{p,j}=\lambda_{p,j}$).

\begin{figure}[t!]
\centering
{\includegraphics[height =6.5cm]{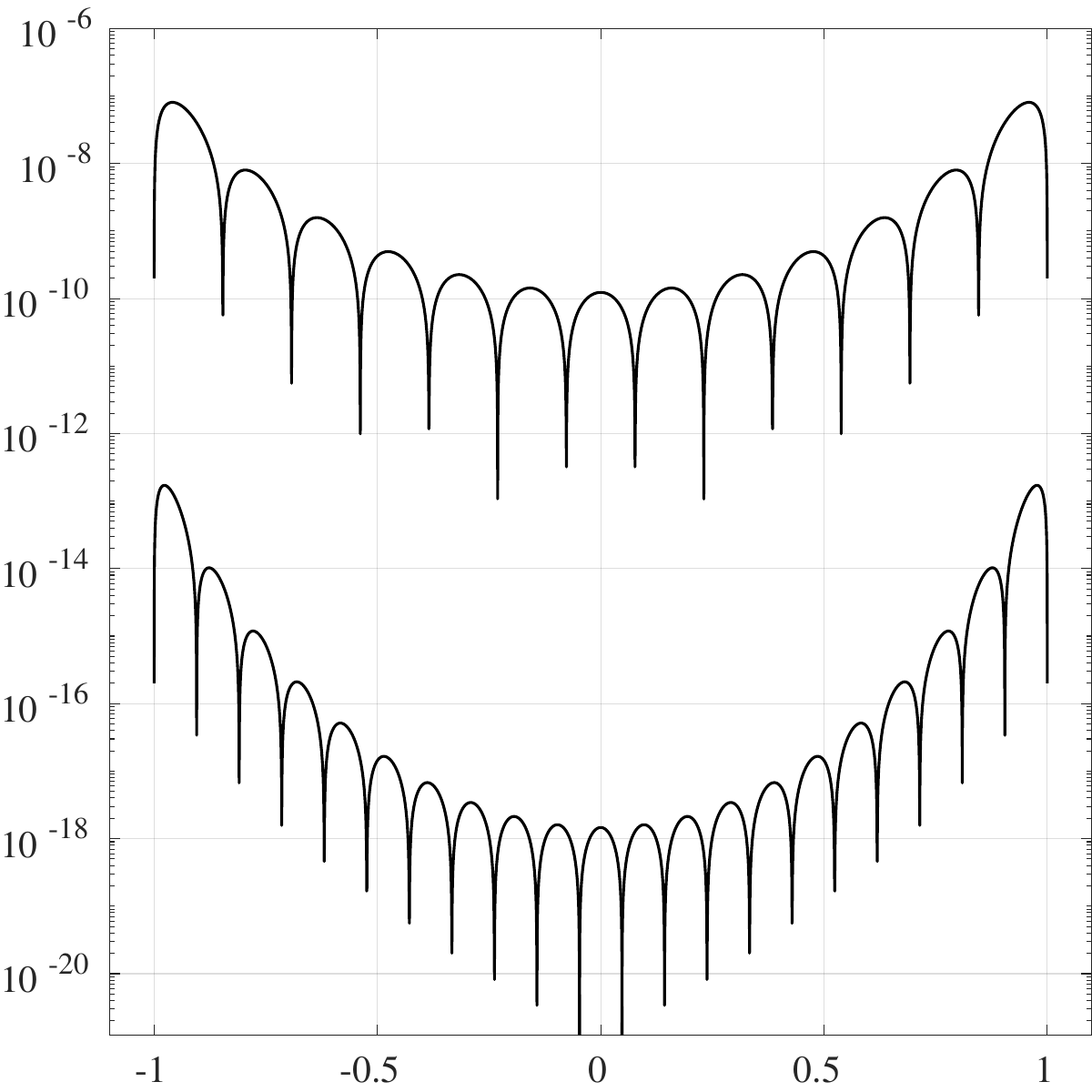}} 
\caption{The values of $|H_p(y) - H(y)|$ for   $-1 \le y \le 1$ and $p \in \{3,5\}$.} 
\label{fig10}
\end{figure}

It remains to explain the pattern in the error 
$\zeta_p(1/2+\i t) - \zeta(1/2+\i t)$ that we observed in
Figure \ref{fig5}.  Define the function
\begin{equation}\label{def:H_p}
H_p(y):=\omega_{p,0} + 2 \sum\limits_{j=1}^p 
\omega_{p,j} e^{-\pi \i \lambda_{p,j}^2} 
\cosh(2\pi \lambda_{p,j} y).
\end{equation}
The system of equations \eqref{eqn:quadrature_condition2} is equivalent to 
\begin{equation}\label{eqn:H_p_H_yk}
H_p(y_{p,k})=H(y_{p,k}), \;\;\; k=0,1,\dots,4p+1,
\end{equation}
where we recall that $y_{p,k}=-1+2k/(4p+1)$ and $H(y)$ is defined in \eqref{eqn:Mordell_integral}. 
We remark here that equations \eqref{eqn:H_p_H_yk} are particularly useful: they allowed us to verify the accuracy of the computed values of $\{\omega_{p,j}\}_{0\le j \le p}$ and $\{\lambda_{p,j}\}_{1\le j \le p}$.
To proceed, we interpret $H_p(y)$ as an approximation to the Mordell integral $H(y)$: this approximation is exact at $4p+2$ equally spaced points on $[-1,1]$. In Figure \ref{fig10}, we plot the error $|H_p(y)-H(y)|$ of this approximation for $p=3$ and $p=5$. These graphs closely resemble those in Figure \ref{fig5}, which display the error $|\zeta(1/2+\i t)-\zeta_p(1/2+\i t)|$ for $p\in \{3,5\}$ at large values of $t$. Let us now explore the reason for this resemblance. 

We rewrite \eqref{def:I_M2} in an equivalent form:
\begin{equation*}
{\mathcal I}_M(s)=M^{-s}  \int_{\r} \frac{e^{-2\pi x^2+2\pi \theta B(t) x}}{\cosh(\pi \theta x)} 
  \d x 
  + M^{-s} \int_{\r} \frac{e^{-2\pi x^2+2\pi \theta B(t) x}}{\cosh(\pi \theta x)}
  \times \Big[e^{\pi x^2- 2\pi \theta B(t) x+g(x;s)} -1 \Big] \d x.
 \end{equation*} 
According to \eqref{eqn:g_asymptotics}, the term inside the square brackets converges to zero as $t \to  +\infty$. 
By the Dominated Convergence Theorem, the second integral on the right-hand side is $o(1)$ as $t\to +\infty$. Thus we obtain 
\begin{equation*}
{\mathcal I}_M(s)=M^{-s} H(B(t)) + o(M^{-s}), \;\;\; t \to +\infty.
\end{equation*}
From \eqref{def:ii_Mp} and \eqref{eqn:g_asymptotics}, we find that as $t\to +\infty$
\begin{align*}
M^{s} \ii_{M,p}(s)&= \omega_{p,0}    +  \sum\limits_{j=1}^p 
\omega_{p,j} \Big[
e^{g(\lambda_{p,j}/\theta;s)} 
 +e^{g(-\lambda_{p,j}/\theta;s)}  \Big] \\&\longrightarrow
 \omega_{p,0}    +  \sum\limits_{j=1}^p 
\omega_{p,j} \Big[
e^{- \pi \i \lambda_{p,j}^2 + 2\pi B(t) \lambda_{p,j} } 
 +e^{- \pi \i \lambda_{p,j}^2 - 2\pi B(t) \lambda_{p,j} }   \Big]
 = H_p(B(t)).
\end{align*}
From these results, it follows that
\begin{equation}\label{I_M_error}
\ii_{M,p}(s)-\ii_{M}(s)=M^{-s} \times \big( H_p(B(t))-H(B(t)) \big) + o(M^{-s}), \;\;\; t\to +\infty,
\end{equation}
which implies
\begin{equation}\label{bar_I_M_error}
\overline \ii_{M,p}(1-s)-\overline \ii_{M}(1-s)=M^{s-1} \times \big( \overline H_p(B(t))-\overline H(B(t)) \big) + o(M^{s-1}), \;\;\; t\to +\infty. 
\end{equation}

Using \eqref{def:chi} and Stirling's asymptotic formula for the Gamma function, we obtain  
$$
\chi(s)=M^{1-2s} \times e^{\frac{\pi \i}{4} + \i t} \Big(\frac{t}{2\pi M^2}\Big)^{\frac{1}{2}-s} \times 
\big(1+O(t^{-1})\big). 
$$
Denote $\tau:=2(\sqrt{t/(2\pi)}-M)$. Writing 
$$
\ln\Big(\frac{t}{2\pi M^2}\Big)=\ln\Big(1+\frac{t-2\pi M^2}{2\pi M^2}\Big)=  \ln\Big(1+\tau \frac{\sqrt{t/(2\pi)}+ M}{2 M^2}\Big), 
$$
expanding the logarithm in Taylor series and simplifying the result, we obtain 
$$
 \chi(s)=M^{1-2s} \times e^{\frac{3 \pi \i}{4}- \pi \i \tau^2}\times \Big(1+
O\big(t^{-\frac{1}{2}}\big)\Big).
$$ 
Combining all the above results with \eqref{def:zeta_p}, 
\eqref{eqn:zeta_as_RS}, \eqref{I_M_error} and  \eqref{bar_I_M_error} we arrive at 
\begin{equation}\label{zeta_p_asymptotic_error}
\zeta_p(s)-\zeta(s)=-\frac{(-1)^N}{2} M^{-s} \Big[ 
 \big( H_p(B(t))-H(B(t)) \big) 
+ e^{\frac{3 \pi \i}{4}- \pi \i \tau^2} \times \big( \overline H_p(B(t))-\overline H(B(t)) \big)\Big]+ o(M^{-s}).
\end{equation}
Thus, we see that the dominant term in the error $\zeta_p(s)-\zeta(s)$ vanishes when $B(t)=y_{p,k}$. This explains the behaviour of the error in Figure \ref{fig5}: when $n$ is large, the error $\zeta_p(s)-\zeta(s)$ is small at $4p+2$ equally spaced points within the interval $[t_n, t_{n+1}]$ because the dominant asymptotic term is zero at those points. 

\section{Comparison with the Riemann-Siegel formula}\label{Section_compare_with_RS}

The Riemann-Siegel approximation to $\zeta(s)$ has long been the preferred way to compute $\zeta(s)$ on and near the critical line $\re(s)=1/2$ when $\im(s)$ is large \cite{Borwein_2000,Gabcke1979,Odlyzko_1988,Rubinstein_2005}. 
In this section, we compare our approximations $\zeta_p(s)$ with the Riemann-Siegel approximations on the critical line. Following \cite{Gabcke1979},
we denote 
\[
\vartheta(t):=\im \big[\ln \Gamma(1/4+\i t/2)\big]-t \ln(\pi)/2,
\]
and for an integer $K\ge 0$ we define the Riemann-Siegel approximation
to $\zeta(1/2+\i t)$ (for $t>0$) as follows: 
\[
\zeta^{\textnormal{RS}}_K(t):= e^{-\i \vartheta(t)}\Bigg [ 2  \sum\limits_{n=1}^{N_t} \frac{\cos(\vartheta(t)-t \ln(n))}{\sqrt{n}}
+\frac{(-1)^{N_t-1}}{\sqrt{a}} \sum\limits_{n=0}^K \frac{C_n(z)}{a^n} \Bigg],
\]
where  
\[
a:=\sqrt{t/(2\pi)}, \;\;\; N_t:=\lfloor a \rfloor, \;\;\; z:=1-2(a-N_t),
\]
and
\begin{align}
\label{def_Fz}
C_0(z)&=F(z):=\frac{\cos(\frac{\pi}{2}(z^2+3/4))}{\cos(\pi z)}, \\
\nonumber
C_1(z)&=\frac{F^{(3)}(z)}{2^2\cdot 3 \pi^2}, \;\;\;
C_2(z)=\frac{F^{(6)}(z)}{2^5\cdot 3^2 \pi^4}+
\frac{F^{(2)}(z)}{2^4 \pi^2}.
\end{align}
All the above formulas, as well as explicit expressions for $C_n(z)$ for $ n \le 12$, can be found in \cite{Gabcke1979}. 
The general form of $C_n(z)$ is  
\begin{equation}\label{formula_Cn}
C_n(z)=2^{-2n} \sum\limits_{k=0}^{\lfloor 3n/4 \rfloor}
\frac{d_k^{(n)}}{\pi^{2n-2k} (3n-4k)!} F^{(3n-4k)}(z),
\end{equation}
where the rational numbers $d_k^{(n)}$ are computed recursively (see formula (3) in \cite{Gabcke1979}). 
It is known that $R_K(t):=\zeta^{\textnormal{RS}}_K(t)-\zeta(1/2+\i t)=O(t^{-(2K+3)/4})$ as $t\to +\infty$ and  Gabcke \cite{Gabcke1979} proves that for $t\ge 200$
\begin{align}\label{RK_bounds}
&|R_0(t)|<0.127 \, t^{-3/4}, \;\; |R_1(t)|<0.053 \, t^{-5/4}, \;\; 
 |R_2(t)|<0.011 \, t^{-7/4},\\ \nonumber
&|R_3(t)|<0.031 \, t^{-9/4}, \;\;\; |R_4(t)|<0.017 \, t^{-11/4},\;\;\; |R_5(t)|<0.061 \, t^{-13/4},
\end{align}
with the bounds for $K\le 4$ being optimal.

\begin{figure}[t!]
\centering
\subfloat[]{\label{fig_11a}
\includegraphics[height =6.35cm]{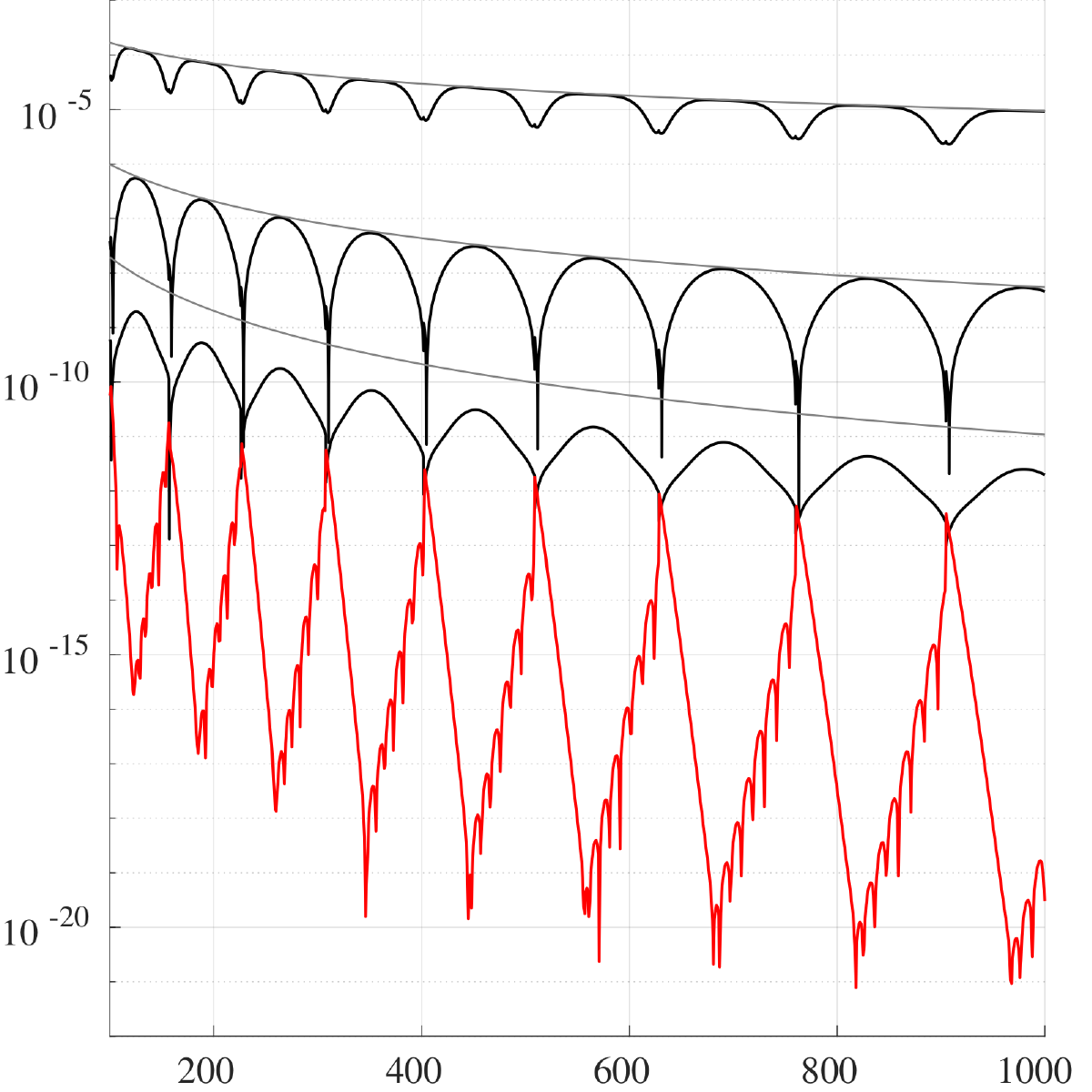}}
\subfloat[]{\label{fig_11b}
\hspace{0.2cm}
\includegraphics[height =6.5cm]{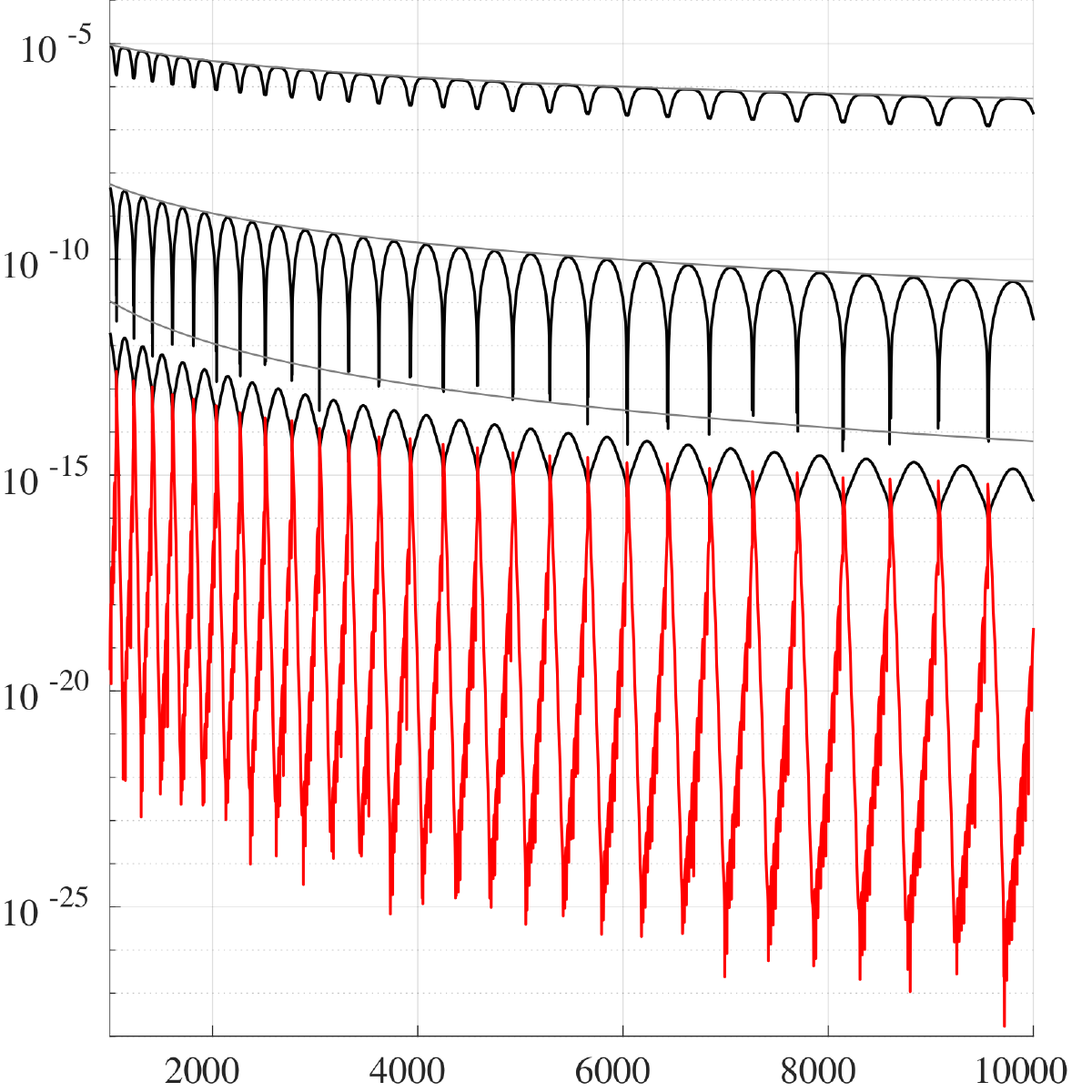}}
\caption{The red curve shows the graph of $|\zeta_7(1/2+\i t)-\zeta(1/2+\i t)|$ ($t$ is on the $x$-axis). The black curves are the Riemann-Siegel approximation errors $|R_K(t)|=
|\zeta^{\textnormal{RS}}_K(t)-\zeta(1/2+\i t)|$ for $K \in \{1,3,5\}$ and the gray curves are the corresponding upper bounds  \eqref{RK_bounds}.} 
\label{fig11}
\end{figure}  

One drawback of the Riemann-Siegel approximations, already noted in the literature (see \cite[Page 259]{Borwein_2000}), is that the correction terms $C_n(z)$ with $n\ge 1$ are not easy to evaluate. When computing $C_n(z)$ we need the derivative $F^{(3n)}(z)$: these can be computed explicitly, but the formulas are very complicated even for small values of $n$. An even more serious problem is that the function $F(z)$ defined in \eqref{def_Fz} has a removable singularity at $z=\pm 1/2$. Note that $F$ is an entire function of $z$ -- the zeros of the denominator $\cos(\pi z)$ are cancelled by the zeros of the numerator $\cos(\pi (z^2+3/4))$. However, the explicit expression for $F^{(j)}(z)$ will have $\cos(\pi z)^j$ in the denominator, which makes it hard to compute the  values of $F^{(j)}(z)$ when $z$ is close to $\pm 1/2$. The standard solution to this problem is not to use the explicit formulas for $C_n(z)$ for $n\ge 1$ and instead to use Taylor or Chebyshev expansions. The coefficients of Taylor and Chebyshev expansions of $C_n(z)$ (given with a precision of 50 decimal digits) can be found in Tables IV and V in \cite{Gabcke1979}.

In Figure \ref{fig11} we plot the errors of the Riemann-Siegel approximations $\zeta^{\textnormal{RS}}_K(t)$ for $K\in \{1,3,5\}$, the corresponding upper bounds \eqref{RK_bounds} (gray curves) and the error of our approximation $\zeta_7(1/2+\i t)$ (red curve). We see that our approximation $\zeta_7$ is better than the Riemann-Siegel approximation $\zeta^{\textnormal{RS}}_5$ in the range $10^2\le t \le 10^3$, and it is better for most values of $t$ in the range $10^3 \le t \le 10^4$. For very large values of $t$, the Riemann-Siegel approximation $\zeta^{\textnormal{RS}}_5$ will be more accurate than $\zeta_7$, as it has a faster rate of asymptotic decay. Indeed, while $R_5(t)$ decays as $O(t^{-13/4})$, formula \eqref{zeta_p_asymptotic_error} implies that the error $|\zeta_7(1/2+\i t)-\zeta(1/2+\i t)|$ decays as $O(t^{-1/4})$. 

At the same time, our approximation $\zeta_7$ is considerably easier to implement than $\zeta^{\textnormal{RS}}_5$. Computing $\zeta_7$ requires only $15$ precomputed complex coefficients $\{\omega_{7,j}\}_{0\le j \le 7}$ and $\{\lambda_{7,j}\}_{1\le j \le 7}$ (these coefficients can be downloaded  \href{https://kuznetsovmath.ca/}{here}). On the other hand,  computing $\zeta^{\textnormal{RS}}_5$ to full double precision requires around 100 coefficients of the Taylor series expansion (20 for each $C_n(z)$ term for $n=1,2,\dots,5$), or  around 65 coefficients of the Chebyshev expansion (13 for each $C_n(z)$ term for $n=1,2,\dots,5$) (see Tables IV and V in \cite{Gabcke1979}). When computing $\zeta^{\textnormal{RS}}_5$ in full quadruple precision, we will need even more: namely, about 170 Taylor or 120  Chebyshev coefficients. Thus, when comparing our approximations $\zeta_p$ with the Riemann-Siegel approximations $\zeta^{\textnormal{RS}}_K$ that achieve similar accuracy, our approximations have the advantage of requiring fewer precomputed coefficients and being easier to implement. The disadvantage is that currently there are no rigorous error bounds for $\zeta_p(s)-\zeta(s)$, but we hope this is a temporary problem and  that efficient and rigorous error bounds will be obtained in the near future.

\section{Computing $\zeta(s)$ in the entire complex plane}\label{Section_computing_zeta}

\begin{figure}[t!]
\centering
\includegraphics[height =6.35cm]{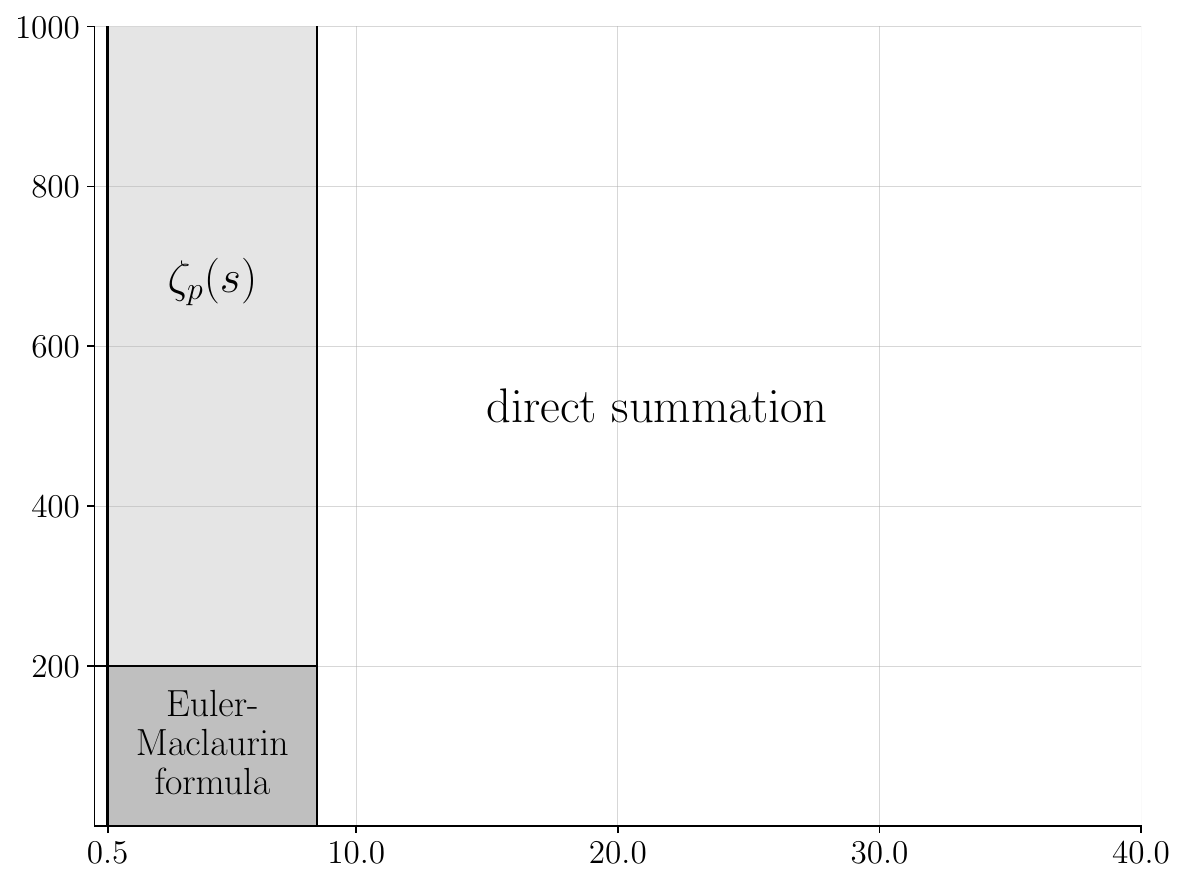}
\caption{The algorithms used to compute $\zeta(s)$ in the region $\im(s)\ge 0$ and $\re(s)\ge 1/2$. } 
\label{fig12}
\end{figure}  

The Riemann zeta function is an important and well-studied special function, so one would expect that most programming languages would have available routines for computing $\zeta(s)$ for arbitrary complex values of $s$. While a number of such routines are available (for example, MATLAB and Maple implementations), they tend to be rather slow and inefficient. In this section, we explain how our approximations $\zeta_p(s)$ can be used together with the direct summation method and Euler-Maclaurin approximations to compute $\zeta(s)$ for arbitrary complex values of $s$.

We note first that we only need to focus on computing  $\zeta(s)$ in the region $\re(s)\ge 1/2$ and $\im(s) \ge 0$, as the functional equation $\zeta(s)=\chi(s) \zeta(1-s)$ and the symmetry relation $\zeta(s)=\overline{\zeta(\overline s)}$ allow one to compute $\zeta(s)$ for every $s\in \c$. In the region
$\re(s)\ge 1/2$ and $\im(s) \ge 0$ we will switch between three different methods for computing $\zeta(s)$, depending on which method is more efficient for the given value of $s$.  

The first method is based on direct summation. We assume that $\re(s)\ge C_1>1$ and write 
\begin{equation}\label{direct_sum}
\zeta(s)=\sum\limits_{n=1}^N n^{-s} + {\mathcal E}_1(N,s),
\end{equation}
where the error ${\mathcal E}_1(N,s)$ is bounded by: 
\[
|{\mathcal E}_1(N,s)|\le \sum\limits_{n>N} n^{-\re(s)}<\int_N^{\infty} x^{-C_1} \d x=\frac{1}{C_1-1} N^{1-C_1}. 
\]
For example, setting $C_1=8$ we find that the direct summation with $N=150$ terms computes $\zeta(s)$ in the half-plane $\re(s)\ge 8$ with error  $\vert {\mathcal E}_1(150,s)\vert < 8.5 \times 10^{-17}$, thus achieving full double precision.

The second method is the well-known Euler-Maclaurin summation formula,  see \cite[Section 3]{Borwein_2000} and \cite[Section 2.2.5]{Rubinstein_2005}, which has the form 
\begin{equation}
\label{zeta_EU}
\zeta(s)=\sum\limits_{n=1}^{N-1} n^{-s}
+\frac{1}{2} N^{-s}
+\frac{N^{1-s}}{s-1}+\sum\limits_{k=1}^K  \frac{B_{2k}}{(2k)!} 
N^{1-s-2k} \prod\limits_{j=0}^{2k-2} (s+j)+ {\mathcal E}_2(s,N,K), 
\end{equation}
As was explained in \cite{Rubinstein_2005}, the error term satisfies $|{\mathcal E}_2(s,N,K)|<10^{-d}$ if $\re(s) \ge 1/2$ and
\[
2\pi N \ge 10 |s+2K-2|, \;\;\; 2K-1 > d+\frac{1}{2} \log_{10} |s+2K-1|.
\]
Thus,
in the region  $1/2 \le \re(s) \le C_1$ and $0\le \im(s) \le C_2$ (this is the dark gray region in Figure \ref{fig12}) we can use the Euler-Maclaurin summation method to compute $\zeta(s)$ to any desired accuracy $\epsilon$ in $O_{\epsilon}(|s|)$ arithmetic operations. 

The third method for computing $\zeta(s)$ is our approximation $\zeta_p(s)$, which will be used in the region $1/2 \le \re(s) \le C_1$ and $\im(s) \ge C_2$ (the light gray region in Figure \ref{fig12}). 

To make the above algorithm efficient, one needs to choose the values of $p$, $C_1$, and $C_2$ carefully. 
For example,  if we take $C_2$ very large, 
then our numerical results in Section \ref{section:Intro} indicate that we can use a smaller value of $p$, so that $\zeta_p(s)$ still gives an approximation to $\zeta(s)$ with the desired accuracy. The downside is that for some values of $s$ with $\im(s)$ large we will have to use the Euler-Maclaurin method, and this will require more computation time. Similar considerations apply to choosing $C_1$: increasing $C_1$ will require fewer terms in the direct summation method but may also require increasing $p$ and require more terms in the Euler-Maclaurin method. Thus, the parameters $p$, $C_1$, and $C_2$
should be chosen to balance the desired accuracy against computational cost.   

We implemented this algorithm for computing $\zeta(s)$ and $\zeta'(s)$ for general complex values of $s$ in MATLAB and Python (in double precision) and in Fortran (in quadruple precision). This code can be downloaded at 
\href{https://github.com/Alexey-Kuznetsov-math/Riemann_zeta-and-Gamma}{GitHub} or from the author's \href{https://kuznetsovmath.ca/}{webpage}. 

For the double-precision MATLAB and Python implementations, we used $p=8$, $C_1=5$ and $C_2=200$. The resulting values of $\zeta(s)$ have typical relative errors of order 
$10^{-13}$ for $|\im(s)|\le 100$, $10^{-12}$ for $|\im(s)|\le 10^3$ and $10^{-11}$ for $|\im(s)|\le 10^4$. 
This loss of precision for large values of $\im(s)$ is consistent with what we observed in the example at the end of Section \ref{section:Intro} (see Figure \ref{fig7}). At the same time, our MATLAB version of $\zeta(s)$ is faster than the MATLAB built-in function ${\texttt{zeta(s)}}$ by a factor of hundreds or thousands  (depending on the range of $|\im(s)|$).

For the Fortran implementations we used $p=30$ and $C_2=400$. We did all calculations in quadruple precision and for every $s$ in the region $0\le \im(s)\le 400$,  we used either Euler-Maclaurin summation or direct summation, depending on which of these has smaller computational complexity for that particular value of $s$. Similarly, in the region $\im(s)>400$ we switched between the $\zeta_{30}(s)$ approximation and direct summation, depending on which has smaller computational complexity. After testing this implementation of $\zeta(s)$, we found that it had typical relative errors of order $10^{-31}$  in the region $|\im(s)| < 100$, $10^{-30}$ for $|\im(s)| < 10^3$, and $10^{-29}$ for $|\im(s)| <10^4$.

\section*{Acknowledgements}

This research was supported by the Natural Sciences and Engineering Research Council of Canada.  The author is
 grateful to two anonymous referees for carefully reading the paper and for providing helpful comments.



\newpage 

\appendix
\section*{Appendix A  \hspace{0.5cm} The coefficients $\omega_{p,j}$ and $\lambda_{p,j}$ for $p \in \{5,10\}$}

High-precision values of coefficients $\{\omega_{p,j}\}_{0\le j \le p}$ and $\{\lambda_{p,j}\}_{1\le j \le p}$ 
for $p=1,2,\dots,30$ and for many values in the range $p\le 150$  can be downloaded from  \url{https://kuznetsovmath.ca/}.

 { \small
  \begingroup
\addtolength{\jot}{-0.2em}

\begin{align*}
\omega_{5,0}&=2.354383173482941501\text{\text{e-}}1+\i \times   3.295537698903362209\text{\text{e-}}2,\\
\omega_{5,1}&=1.737747054311600606\text{e-}1+\i \times   5.726284240637533629\text{e-}2, \\
\omega_{5,2}&=5.708483151712981717\text{e-}2+\i \times   5.367826163392843596\text{e-}2, \\
\omega_{5,3}&=5.455260017239278090\text{e-}3+\i \times    1.692893836426674071\text{e-}2, \\ 
\omega_{5,4}&=-4.138943380524367182\text{e-}4+\i \times    2.034356935140766843\text{e-}3, \\
\omega_{5,5}&=-6.653714335968984044\text{e-}5+\i \times    6.462574635932469471\text{e-}5       
    \end{align*}            
\begin{align*}
\lambda_{5,1}&=1.881852180702220422\text{e-}1-\i \times    1.449755745395875068\text{e-}1 , \\
\lambda_{5,2}&=3.717705006684543749\text{e-}1 -\i \times    3.020416160270775712\text{e-}1, \\
\lambda_{5,3}&=5.604567753443639984\text{e-}1-\i \times    4.801938441334826040\text{e-}1 , \\ 
\lambda_{5,4}&=7.672689586414600920\text{e-}1-\i \times     6.821905910209796031\text{e-}1 , \\
\lambda_{5,5}&= 1.010783983564685329-\i \times     9.231734623461863512\text{e-}1.       
    \end{align*}

\begin{align*}
\omega_{10,0}&=1.746071737157674980979293520809\text{e-}1 +\i \times    					   2.131147093009280730611467019158\text{e-}2  ,\\
\omega_{10,1}&=1.490803915553910597329354639778\text{e-}1+\i \times    					   3.499836079601156948133789078972\text{e-}2, \\
\omega_{10,2}&=8.492465921092508217336004148263\text{e-}2+\i \times    					   4.854991766416009886502556092917\text{e-}2, \\
\omega_{10,3}&=2.794492162555768303150174880103\text{e-}2+\i \times    					   3.428439466181300925395192520791\text{e-}2  , \\ 
\omega_{10,4}&=4.612090699061725829646273271703\text{e-}3+\i \times      					   1.373142646307427391022925045066\text{e-}2, \\
\omega_{10,5}&=-3.895212927973588318860893961158\text{e-}5+\i \times      			    3.550886924259579942806268192521\text{e-}3 , \\
\omega_{10,6}&=-2.151575611923250640729364801406\text{e-}4+\i \times      				6.084356024918800989143391852680\text{e-}4 , \\
\omega_{10,7}&=-5.199488450834904743451274940186\text{e-}5+\i \times      			    6.406830664562431793000193930144\text{e-}5, \\
\omega_{10,8}&=-5.856003331642731075366848061989\text{e-}6+\i \times      			    3.353733365341979352981823198386\text{e-}6, \\ 
\omega_{10,9}&=-2.945578758160111306783176275407\text{e-}7+\i \times      			    3.154278990732981364449273807939\text{e-}8, \\
\omega_{10,10}&=-4.219551146037265608639695765718\text{e-}9-\i \times      				 1.752142489214440816303376939714\text{e-}9,        
    \end{align*}            
\begin{align*}
\lambda_{10,1}&= 1.379409313309054508271675868217\text{e-}1  -\i \times      				 1.088692797924869220391271752962\text{e-}1 , \\
\lambda_{10,2}&= 2.732463550335757861584970430657\text{e-}1      -\i \times           2.210503737259508831029856904771\text{e-}1 , \\
\lambda_{10,3}&= 4.070334053056538299722767959949\text{e-}1  -\i \times      				 3.400869979247635282520012627532\text{e-}1  , \\ 
\lambda_{10,4}&= 5.429713841237013800653833464349\text{e-}1   -\i \times      				 4.668200118355472525024744280421\text{e-}1  , \\
\lambda_{10,5}&= 6.834620082884849199273619613380\text{e-}1 -\i \times      				 6.002854340275813175341293481950\text{e-}1  ,\\  
\lambda_{10,6}&= 8.297493681957483741659306681846\text{e-}1 -\i \times      				 7.404377940784227473659159034325\text{e-}1 , \\
\lambda_{10,7}&= 9.835018784062355446404273245147\text{e-}1  -\i \times      				 8.888012731779453622778359704903\text{e-}1 , \\
\lambda_{10,8}&= 1.147933282145432947538394279481   -\i \times      						     1.048670473139049661794532732170, \\ 
\lambda_{10,9}&= 1.329633190044527778848402344442    -\i \times      							 1.226639730249438411182742778670, \\
\lambda_{10,10}&=1.545989175497797759478691005072     -\i \times      					     1.440017038829556195286509733096.      
    \end{align*}     
   
      \endgroup 
      }
      
     \newpage 
\appendix
\section*{Appendix B  \hspace{0.5cm} MATLAB code for $\zeta_8(s)$}\label{AppendixB}
 This function was tested  numerically for $1/2\le \re(s)\le 2$ and $100\le \im(s)\le 10000$, where it computes values of $\zeta(s)$ with an accuracy of $10$ to $12$ decimal digits. The accuracy decreases as $\im(s)$ increases, primarily due to rounding errors (see Figure \ref{fig7}). Implementing this code in quadruple precision solves the problem with rounding errors (unless $t$ is very large, of order $10^{15}$ or greater).  
   \begin{lstlisting}[
frame=single,
numbers=left,
style=Matlab-editor,
basicstyle=\footnotesize]
function f=zeta_8(s)
% this function computes zeta_8(s) for imag(s)>0
% the cofficients lambda_{8,j} for j=1,2,...,8
lambda=[0.152845417613666702426-0.119440685603870510384i
        0.302346225128945757427-0.243989695504400621268i          
	0.451119584531782942888-0.378479770209444563858i 
	0.604563710297226464637-0.523486888629095259770i
	0.765965706759629396959-0.678405572413543444272i
	0.938371150977889047740-0.845332361280975174880i
	1.128148837845288402558-1.030737947568157685685i
	1.353030558654668162533-1.252503278108132307164i];    
% the cofficients omega_{8,j} for j=0,1,2,...,8	           
omega0=1.926019633029103199063e-1+2.472986965795651842299e-2i;
omega=[1.582954327321094104502e-1+4.149113569204600502105e-2i   
       7.826728293587305110862e-2+5.215518667623989653254e-2i   
       1.940595049247490540621e-2+2.977286598777633378610e-2i  
       1.691184771902755036966e-3+8.938933548999206800196e-3i     
      -2.994777986686168319731e-4+1.567541981830224487301e-3i   
      -9.837202592542590210980e-5+1.502108057352792742070e-4i  
      -9.346989286415688998740e-6+5.793852209955845432028e-6i  
      -2.451577304299235983015e-7+6.134784898751456953524e-9i];  
% compute chi(s)=(2*pi)^s/(2*cos(pi*s/2)*gamma(s))
% we use Stirling's formula for log(gamma(s)) and truncate
% the asymptotic series \sum_{n\ge 1} B_{2n}/(2*n*(2*n-1)*s^(2*n-1)) at n=3
chi=exp((0.5-s)*log(s/(2*pi))+s*(0.5i*pi+1)-(1/s)*(1/12+s^(-2)*(-1/360+s^(-2)/1260)));
% compute the main sum
N=floor(sqrt(imag(s)/(2*pi)));
lnn=log(2:N); 
f=1+sum(exp(-s*lnn))+chi*(1+sum(exp((s-1)*lnn)));   
% compute I1=I_{M,8}(s) 
M=N+0.5;
I1=exp(-s*log(M))*(omega0+sum(omega.*(exp(-2*pi*M*lambda-s*log(1+1i*lambda/M))          +exp(2*pi*M*lambda-s*log(1-1i*lambda/M)))));  
% compute I2=\bar I_{M,8}(1-s)=conj(I_{M,8}(1-conj(s)))
z=1-conj(s); 
I2=conj(exp(-z*log(M))*(omega0+sum(omega.*(exp(-2*pi*M*lambda-z*log(1+1i*lambda/M))      +exp(2*pi*M*lambda-z*log(1-1i*lambda/M))))));
% compute zeta_8(s)
f=f-0.5*(-1)^N*(I1+chi*I2); 



   \end{lstlisting}

\end{document}